\documentclass[a4paper,12pt,final]{article}
    \usepackage{tikz}
    \usepackage{xcolor}
    \usepackage[utf8x]{inputenc}
    \usepackage[T1]{fontenc}
    \usepackage[scr=boondoxo,scrscaled=1.05]{mathalfa}
    \usepackage{lmodern}
    \usepackage{textcomp}
    \usepackage[inkscapearea=page]{svg}
    \usepackage{adjustbox}
    \usepackage[a4paper,pdftex,dvips,left=2.2cm, right=2.2cm, top=2.7cm, bottom=2.7cm]{geometry}
\usepackage{amsmath,amsthm,amsfonts,amssymb,amscd,url,mathtools}
    \usepackage[backref=page,colorlinks=false]{hyperref}
    \usepackage{epsfig} 
    \usepackage{faktor}
    \usepackage{xfrac}
    \usepackage{float}
    \usepackage[english,french]{babel}
    \usepackage{graphics}
    \input{xy} \xyoption{all}
    
    \usepackage{enumerate,xcolor}
    \usepackage{mathrsfs}
    \usepackage{epsfig}
    \usepackage[normalem]{ulem}
    \usepackage[toc,page]{appendix}
    \usepackage{pgfplots}
    \usepackage[scr=boondoxo,scrscaled=1.05]{mathalfa}
    \makeindex

    \def\diam {\mathrm{diam\, }}
    \def\radius {\mathrm{radius\, }}

    \def\exp{\mathrm{exp}}

    \def\loc{\mathrm{loc}}

    \def\R{\mathbb R}

    \def\E{\mathbb E}

    \def\diam{\mathrm{diam}\:}

    \definecolor{blue}{rgb}{0,0,225}

    \def \exp{\mathrm{exp}}

    \def\qand{\quad \text{and}\quad}

    \def\Q{\mathbb Q}

    \def\R{\mathbb R}

    \def\diam{\mathrm{diam}\:}

    \def\qand{\quad \text{and} \quad}
    
    \def\1{{1\!\! 1}}
    \theoremstyle{plain}

    \usepackage{array, babel, booktabs, cite, float, footmisc, kvoptions,multicol, nag, newtxmath, newtxtext, pdf14, pdftexcmds, ragged2e, upref, url, xcolor, xpatch, zref-base}
    \usepackage[capitalize]{cleveref}
\theoremstyle{plain}
\newtheorem{thm}{Theorem}[section]
\newtheorem{theorem}[thm]{Theorem}
\newtheorem{question}[thm]{Question}
\newtheorem{proposition}[thm]{Proposition}

\newtheorem{corollary}[thm]{Corollary}
\newtheorem{lemma}[thm]{Lemma}
\newtheorem{remark}[thm]{Remark}
\newtheorem{fact}[thm]{Fact}

\newtheorem{definition}[thm]{Definition}
\newtheorem{theo}{Theorem} 

\crefname{theo}{Theorem}{Theorems}
\crefname{lemma}{Lemma}{Lemmas}
\crefname{prop}{Proposition}{Propositions}
\crefname{corollary}{Corollary}{Corollaries}
\crefname{theorem}{Theorem}{Theorems}
\crefname{definition}{Definition}{Definitions}
\crefname{remark}{Remark}{Remarks}
\crefname{example}{Example}{Examples}
\crefname{conjecture}{Conjecture}{Conjectures}
\crefname{claim}{Claim}{Claims}
\crefname{problem}{Problem}{Problems}
\crefname{fact}{Fact}{Facts} 
\crefname{quest}{Question}{Questions}
\crefname{conj}{Conjecture}{Conjectures}

    \renewcommand*{\backref}[1]{}
    \renewcommand*{\backrefalt}[4]{\quad \tiny 
      \ifcase #1 (\textbf{NOT CITED.})%
      \or    (Cited on page~#2.)%
      \else   (Cited on pages~#2.)%
      \fi}
    \usepackage{libertine}
    \newcommand{\supess}{ \mathrm{ess \ sup\,  }}
    \newcommand{\infess}{ \mathrm{ess \ inf\,  }}
    
    \def \Scl{\mathsf{scl}}
    \def \scl{\mathit{scl}}

    \def \diam{\mathrm{diam}}
    \def \Var{\mathrm{Var }}
    
    \def \dim{\mathsf{dim}}
    \def \ord{\mathsf{ord}}

    \def\crit{\mathrm{crit}}
    \def \R{\mathbb{R}}
    
    \usepackage{etoolbox}
    \pgfplotsset{compat=1.18} 
    
\begin{document}
    \selectlanguage{english}
    \title{Infinite Dimensional Multifractal Analysis of the Wiener measure}
    \author{Aihua Fan
    \footnote{LAMFA, UMR 7352, University of Picardie, 33 Rue Saint Leu, 80039, Amiens, France and Institute for Math $\&$
AI, Wuhan, Wuhan University, Hubei, 430072, China. \ {\em Email address}: ai-hua.fan@u-picardie.fr}
    \& Mathieu Helfter \footnote{
Institute of Science and Technology Austria (ISTA), Klosterneuburg, Lower Austria. Partially supported by ERC SPERIG $\#885707$. \ {\em Email address}: mathieu.helfter@ist.ac.at
}
}

\date{}
\maketitle
    \selectlanguage{english}
    \begin{abstract}
We present a multifractal formalism for measures on infinite 
dimensional metric spaces, in terms of scales instead of dimensions in 
the classical multifractal analysis. We prove a multifractal formalism 
with a suitable scaling, called order, for the Wiener measure, which 
is the probability law of the standard Brownian motion. We also prove 
the fundamental Frostman Lemma on a large class of Polish spaces, for which the increasing sets lemma  holds.
    \end{abstract}

{ \small  \tableofcontents}

    \section{Introduction}
    
    Multifractal analysis is an important tool in understanding the fine structure of measures, providing a detailed description of local behaviors. A rich theory has been developed in the finite-dimensional setting, most notably for measures on the Euclidean space $ \R^n$ where multifractal spectra are obtained for a wide variety of deterministic and stochastic special measures. Classical examples include self-similar measures \cite{falconer1994,olsen1995}, Gibbs measures associated with dynamical systems \cite{pesin1997}, and stochastic processes such as Bernoulli convolutions \cite{peres1996} and branching processes \cite{barral2002}.
    Given the importance of infinite-dimensional settings in probability theory, analysis, and mathematical physics, we propose to bring  multifractal analysis to a wider landscape describing infinite and zero dimensional measures at their adapted \emph{scales} as introduced in \cite{helfter2025scales}. Scales replace dimensions like Hausdorff dimension and Tricot dimension (also called packing dimension). 
   
Our work is twofold. First, we provide a consistent framework under which the multifractal formalism extends to any scale, and we  establish a Frostman Lemma for arbitrary gauge functions
on a class of Polish spaces including $\sigma$-compact metric spaces, which in turn yields a mass distribution principle valid for any scale. Then, we investigate the first paradigmatic example of infinite-dimensional measures exhibiting multifractal structures, namely the Wiener measure describing standard Brownian motion. 
While several works,  notably   \cite{hua1997multifractal,mannersalo2010multifractal,denisov2016limit}, describe the dimensional multifractals in infinite
product spaces, these often retain a finite  mean dimension and can still be treated within a finite-dimensional formalism. In contrast, the Wiener measure lives on function spaces of genuinely infinite  dimension,  for which classical multifractal spectra degenerate, and new tools are required to capture local fluctuations. We propose to use scales as such tools.
    Brownian motion is a natural choice
    for a simple reason. The space of continuous functions is a topological group, but it fails to be locally compact, which in particular prevents the Wiener measure from playing the role of a Haar measure. This observation already indicates that the Wiener measure would differently concentrate on a large collection of non-typical trajectories, whose irregularities give rise to its multifractal nature. 
   Our investigation will confirm this.
    Our work in the present paper suggests that multifractal analysis in infinite dimensional metric spaces may offer a fruitful perspective on the geometry of stochastic processes and functional spaces.

   Our main results will be collected in Theorems A-G, which serve as the core statements of the paper, while intermediary statements and corollaries use numerical labels.
 \smallskip
In the rest of this introduction, we will
present in some detail, part of these results, namely Theorem A-E, with a brief presentation of the involved notions. In Section 2 we recall the notion of scale and develop the multifractal formalism at any scale, thus generalizing the dimensional case. In Section 3 the generalized Frostman Lemma is proved, and a mass distribution principle is obtained as a corollary which is then applied to study the size of level sets and multifractal spectra. Then in Section 4 we perform the
multifractal analysis of the Wiener measure. This analysis involves functional analysis, geometric
measure theory and stochastic analysis.
Finally in Section \ref{sec:quest} we discuss some questions that naturally  arise from the previous study on the Wiener measure. 
   
\subsection{A general Frostman Lemma}
Frostman Lemma is  fundamental in geometric measure theory. In the case of  Euclidean space $\R^n$ (see \cite{frostman1935potentiel}),  it states that a compact set $K$ has positive $\alpha$-dimensional Hausdorff measure $(\alpha >0)$ if and only if $K$ supports a probability  measure $\mu$  satisfying $\mu(B(x,r)) \leq C r^\alpha$ for balls $B(x,r)$ of radius $r$ and center $x$ ($C>0$ being a constant).
The ``if part'' is trivial but very useful.
Since the Hausdorff measure on $\mathbb{R}^n$ is regular, the result holds not only for compact sets but also for Borel sets, even for analytic sets.
It has been a problem to generalize Frostman Lemma to other metric spaces and to other gauge functions than $r\mapsto r^\alpha$. R. Kaufman \cite{kaufman1994min} proved a weaker version of Frostman Lemma for a compact set $K$ in an arbitrary metric space: 
$\mu(B(x,r)) \leq C r^\alpha \log^{1+\delta}\frac{1}{r}$.
We will improve Kaufman's result by dropping the factor $\log^{1+\delta} \frac{1}{r}$.  Actually, we will prove Frostman’s Lemma for arbitrary gauge functions.

Gauge functions are also called Hausdorff functions. 
We adopt the following definition (see \cite{mattila1997measure},\cite{rogers1998hausdorff}). 

\begin{definition}[Hausdorff function]
A \emph{Hausdorff function}, or \emph{gauge}, is a continuous non decreasing function $ h : [ 0, +\infty)  \to [ 0, +\infty) $ such that $h ( 0) = 0 $ and $ h > 0 $  on $(0,+\infty) $. 
\end{definition}
Given a Hausdorff function $ h $ we will denote $\mathcal{H}^h$ its corresponding $h$-Hausdorff measure. For its definition,  see for instance \cite{rogers1998hausdorff,falconer1997techniques,falconer2004fractal,mattila1997measure}.
We will also recall the definition in Section   \ref{lieu:haus-pack}.

Our first main result establishes a version of Frostman Lemma which is valid for any Hausdorff function and for a large class of separable metric spaces.  

\begin{theo}[Frostman Lemma] \label{thm:frostman-intro}
    Let $ (X,d) $ be a $\sigma$-compact metric space and $ A \subset X$ a Borel subset. Then for every Hausdorff function $h$  such that $\mathcal{H}^h(A) > 0 $, there exist a compact set $K\subset A$,
a  probability Borel measure $\mu$ supported by $K$, and a constant $C\in(0,\infty)$
such that 
\begin{equation*} 
\mu\bigl(B(x,r)\bigr)\ \le\ C\,h( 6 r) , \qquad\forall\,x\in X,\ \forall\,r>0.
\end{equation*}  
\end{theo}
 Actually, we will prove a stronger result (Theorem \ref{thm:Frostman-gauge}), which  states that Frostman Lemma actually remains valid under the mild topological assumption that  the \emph{increasing sets lemma} (see \eqref{eq:increasing-set-lemma}) holds for the Hausdorff measure. This assumption holds for every Hausdorff function if the considered space is, for instance, $\sigma$-compact, according to \cite{sion1962approximation}. Theorem \ref{thm:frostman-intro} is then a direct corollary of this more general statement. 
Theorem \ref{thm:frostman-intro} provides the foundation  for extending  the mass distribution principle (also known as Billingsley theorem) and enables multifractal formalism  beyond the  scalings of power-like gauges.

\subsection{Scalings and scales}
We introduce here a \emph{multifractal formalism} in infinite dimensional metric spaces using \emph{scales} instead of dimensions. We first briefly recall framework of \emph{scales} as introduced in \cite{helfter2025scales}, which are generalizations of dimensions.

A \emph{scaling} is, by definition,  a family $  \Scl =(\mathrm{scl}_\alpha)_{\alpha>0}$ 
of Hausdorff functions satisfying a mild decreasing condition relative to the index $ \alpha $ (see Definition \ref{def:scaling}). 
Typical examples are:
\[
\dim :=(\varepsilon\mapsto \varepsilon^{\alpha})_{\alpha>0}
\qand 
\ord :=(\varepsilon\mapsto \exp (-\varepsilon^{-\alpha}) )_{\alpha>0} \; . 
\] 
Given a scaling $ \Scl = ( \scl_\alpha)_{ \alpha >0 } $ and a metric space $(X,d) $, we denote by $ \mathcal{H}^{\scl_\alpha} $ and $ \mathcal{P}^{\scl_\alpha} $ 
the Hausdorff and packing measures associated to the usual Hausdorff function $ \scl_\alpha $ (see Section \ref{lieu:haus-pack}). For any set $A \subset X $ we define its Hausdorff scale and packing scale as: 
\[ \Scl_H A := \sup \left\{ \alpha > 0 : \mathcal{P}^\scl_{\alpha} (A) = + \infty   \right\} = \inf \left\{ \alpha > 0 : \mathcal{H}^\scl_{\alpha} (A) = 0  \right\} \] 
and
\[ \Scl_P A := \sup \left\{ \alpha > 0 : \mathcal{P}^\scl_{\alpha} (A) = + \infty  \right\}  = \inf \left\{ \alpha > 0 : \mathcal{P}^\scl_{\alpha} (A) = 0  \right\}  \;.  \]
The box scales were also defined (see \cite{helfter2025scales}).

Now, consider a finite Borel measure $\mu$ on 
$(X,d)$ and a scaling $\Scl$,  we define the \emph{local lower and upper scales} of $\mu$ at a point $x \in X $ respectively as 
\begin{equation*}
\underline{\Scl}_{\loc}\mu(x)
   := \sup \left\{\,\alpha>0:\,
     \frac{\mu(B(x,\varepsilon))}{\scl_\alpha(\varepsilon)} \,\underset{\varepsilon\to 0}{\to} 0
  \right\} 
    \end{equation*}
and
\begin{equation*}
  \overline{\Scl}_{\loc}\mu(x)
   := \inf \left\{\,\alpha>0:\,
     \frac{\mu(B(x,\varepsilon))}{\scl_\alpha(\varepsilon)} \,\underset{\varepsilon\to 0}{\to} + \infty 
  \right\} \;. 
\end{equation*}
Whenever they share the same value we denote $ \Scl_\loc \mu (x)$ 
their common value, called the local scale at $x$. 

The case of $ \Scl = \dim $ gives the well known classical Hausdorff, packing and local dimensions. 
The case of  $\Scl =\ord $ gives the so-called ``order''. The local scale describes the small-ball behavior of the measure.  Hausdorff and packing orders were new geometric invariants introduced in \cite{helfter2025scales}. 
 As we are going to see, the order is the suitable scaling for the study of Wiener measure. 

\medskip
When $\Scl_{\loc}\mu(x) $ is $\mu$-everywhere constant, the common value equals to the Hausdorff and packing scales of $\mu$  which must coincide. A more detailed description is recalled below in Theorem \ref{thm:b-scales}.

For a Borel subset $ A \subset X $ let us denote $ \mathcal{M}^1 ( A) $ the set of Borel probability measures supported on $A$. As a consequence of Frostman Lemma (Theorem \ref{thm:frostman-intro}) we have the following Mass distribution principle. 
\begin{theo}[Mass distribution principle]
\label{thm:mass-distribution-scale-intro}
Let $ 
 \Scl = ( \scl_\alpha)_{ \alpha > 0 } $ be a scaling.
 Let $(X,d)$ be a $\sigma$-compact metric space. 
Then for any Borel  $ A \subseteq X$ the Hausdorff scale verifies the mass distribution principle
\[
\Scl_H A = \sup_{\mu \in \mathcal{M}^1(A)}  \inf_{ x \in A} \underline{\Scl}_\loc \mu  (x) \;. 
\]
\end{theo}

\subsection{Multifractal formalism}

The objects of  multifractal analysis are the \emph{level sets} of the local scales of a measure $\mu$ and the objective is to measure their sizes  using Hausdorff or packing scales; thus defining the corresponding \emph{multifractal spectrum}.

    Let us first recall the classical setting of multifractal analysis in finite dimension spaces. Let \( \mu \) be a Borel probability measure on \( \mathbb{R}^d \) with $d \ge 1$. 
     The yet classical multifractal analysis focuses on the description of the distribution of the \emph{local dimensions} of \( \mu \), defined as:
    \[ \underline{\dim}_{\loc} \mu(x) := \liminf_{r \to 0} \frac{\log \mu(B(x, r))}{\log r} \qand \overline{\dim}_{\loc} \mu(x) := \limsup_{ r \to 0} \frac{\log \mu(B(x, r))}{\log r} \;. 
    \]
    When these quantities coincide, we denote by $ \dim_{ \loc} \mu (x) $ their common value. Let:
    \[
    \underline{E} (\alpha) := \left\{ x \in \R^d \;:\; \underline{\dim}_{\loc}\mu(x) = \alpha \right\}  \ ,  \quad \overline{E} (\alpha) := \left\{ x \in \R^d \;:\; \overline{\dim}_{\loc}\mu(x) = \alpha \right\} 
    \]
    and  
    \[
    E(\alpha) :=  \underline{E} (\alpha) \cap \overline{E} (\alpha) \; . 
    \] 
 The main objective of study in this framework is to measure the sizes of these level sets using different dimensions. What we obtain is called the  \emph{multifractal spectra} of \( \mu \), for example, the following  spectra defined by
    \[ \underline{f} ( \alpha) := \dim_H \underline{E} ( \alpha)  \ , \quad 
    \overline{f} (\alpha) := \dim_P \overline{E}(\alpha)  \qand f ( \alpha) := \dim_H E ( \alpha)
    \]
    where \( \dim_H \) denotes the Hausdorff dimension and  \( \dim_P \) denotes the packing dimension.  When $ E ( \alpha) = \emptyset $, the convention  $f ( \alpha) = - \infty $ is often used to distinguish  non-empty sets of $0$ dimension from empty sets.

Now we propose a multifractal formalism for studying measures on an arbitrary metric space by using suitable scalings.
Given a  scaling $ \Scl $ and a real positive number $ \alpha > 0$, we consider: 
\[ \underline{E}^\Scl_\mu (\alpha )  := \{x \in X: \underline{\Scl}_{\loc}\mu(x) =\alpha\},  \quad \overline{E}^\Scl_\mu (\alpha )  := \{x \in X: \overline{\Scl}_{\loc}\mu(x) =\alpha\} \] 
and
\[ E^\Scl_\mu (\alpha ) :=  \underline{E}^\Scl_\mu (\alpha )  \cap \overline{E}^\Scl_\mu (\alpha )    \;.  \]
Even when the local scaling $\Scl_{\loc}\mu$ is well defined and takes a single value $\mu$-almost everywhere, the exceptional sets $E^\Scl_\mu(\alpha)$ for other values of $\alpha$ need not be empty, and may in fact have large Hausdorff or packing scale. These sets of atypical local behavior already play a fundamental role in classical finite-dimensional multifractal analysis and their size is described by the so-called (fine) multifractal spectra.

We thus extend the definition of \emph{multifractal spectra} at a given scaling $ \Scl$ by: 
\[ \underline{f}_\mu^\Scl (\alpha) := \Scl_H \underline{E}_\mu^\Scl ( \alpha)  \, , \quad \overline{f}_\mu^\Scl (\alpha) := \Scl_P \overline{E}_\mu^\Scl ( \alpha)  \qand f_\mu^\Scl (\alpha) = \Scl_H  E_\mu(\alpha)  \; 
\]
if the corresponding level sets are non empty and $ - \infty $ otherwise. 
Note that when $\Scl = \dim $ we indeed retrieve the usual notions of the so-called level sets and (fine) multifractal spectra. 

As a consequence of the mass distribution principle (Theorem \ref{thm:mass-distribution-scale-intro}) and a result from \cite{helfter2025scales} (see also Theorem \ref{thm:b-scales}), we get the following principle for multifractal formalism.
\begin{corollary} \label{coro:mass-distrib-pple-intro}
Let $ \mu$ be a Borel measure on   metric space $ (X,d) $ and let $ \Scl $  be a scaling. Then for every $ \alpha > 0$, the multifractal spectra for $ \Scl $ at $ \alpha > 0 $ verify
$$  \underline{f}^\Scl_\mu ( \alpha)    \ge    \sup   \left\{    \Scl_H \nu  \, : \, \nu \in \mathcal{M}^1 \left( \underline{E}^\Scl_\mu(\alpha) \right) \right\} \, ,  $$ $$ \overline{f}^\Scl_\mu ( \alpha)    \ge      \sup   \left\{    \Scl_P \nu  \, : \, \nu \in \mathcal{M}^1 \left( \overline{E}^\Scl_\mu(\alpha) \right) \right\} $$
and
$$    f^\Scl_\mu ( \alpha)    \ge  \sup   \left\{    \Scl_H \nu  \, : \, \nu \in \mathcal{M}^1 \left( E^\Scl_\mu(\alpha) \right) \right\} . $$
If, moreover, $X$ is $ \sigma$-compact then we have the equalities: 
$$  \underline{f}^\Scl_\mu ( \alpha)    =     \sup   \left\{    \Scl_H \nu  \, : \, \nu \in \mathcal{M}^1 \left( \underline{E}^\Scl_\mu(\alpha) \right) \right\}   $$ and $$ f^\Scl_\mu ( \alpha)   =   \sup   \left\{    \Scl_H \nu  \, : \, \nu \in \mathcal{M}^1 \left( E^\Scl_\mu(\alpha) \right) \right\} . $$
\end{corollary}
Here $\Scl_H\nu$ denotes the Hausdorff scale of the measure $\nu$ (see \eqref{Def:sclH}  for its definition). 
For this principle to be effective, we can try to find a measure $\nu$ which maximizes the above equalities.
For the study the multifractal spectrum of the Wiener measure, we will find that  the probability laws of fractional Brownian motions are useful. In the well known works on the multifractal analysis of Birkhoff averages in dynamical systems, Gibbs measures play  this role (see \cite{fan1994dimensions}, \cite{Fan2021}, \cite{pesin2001multifractal,olsen2003multifractal}). 
 Theorem \ref{thm:mass-distribution-scale-intro}, as well as Corollary \ref{coro:mass-distrib-pple-intro}, relies on the generalization of Frostman  lemma.  We will  provide, for both, stronger versions of the results, see Theorem \ref{thm:mass-distribution-scale} and Corollary \ref{coro:mass-distrib-pple}. 
\subsection{Multifractal formalism of the Wiener measure}

We propose the Wiener measure  as a first paradigmatic example in the study of multifractal analysis in the infinite-dimensional context. The Wiener measure is defined  on the space of continuous functions defined on the interval $[0,1]$, equipped with the  uniform norm.
Let $  B := ( B_t)_{ t \in [ 0,1]}$ denote the standard Brownian motion on the interval $ [ 0, 1 ] $  and let $ W $ denote  the corresponding Wiener measure on the space $ \Omega $ of continuous functions on this interval, taking the value $0$ at $ 0$. For a given norm $ \| \cdot \|$ on the path space $ \Omega$, a small ball probability  behavior is usually defined as the asymptotical behavior of: 
\[ 
W  \left(   B_{ \| \cdot \| } ( \bar{\omega} , \varepsilon) \right)\quad \text{ as} \, \varepsilon \to 0 , 
\] 
where $  B_{ \| \cdot \| } ( \bar{\omega} , \varepsilon) $ denotes the open ball of radius $ \varepsilon $ centered at $ \bar{\omega}\in \Omega $ for the norm $ \| \cdot \| $.

Small ball probabilities for Gaussian (or more general) stochastic processes form a well-developed area of probability theory. In the Gaussian case, the problem is typically to determine the asymptotic behaviour of small ball probability which is closely tied to the geometry of the reproducing kernel Hilbert space (RKHS) of the process and to entropy or eigenvalue estimates for the embedding of this RKHS into the ambient normed space.  
For Brownian motion and related Gaussian processes, classical results (see for instance \cite{shao1993note,talagrand1994small,kuelbs1995small,li1999approximation,LiShao2001,dereich2003link,lifshits1995}) provide upper and lower bounds on these probabilities in a wide range of norms, including the uniform norm, Hölder norms, and $L^p$ norms.  
This is exactly the object of study of local scales.  For a large class of metrics, these asymptotics are studied using (implicitely) the specific scaling of order: 
\begin{equation} \label{def:ord}
\ord  = ( \varepsilon \mapsto \exp( - \varepsilon^{-\alpha} ) )_{ \alpha > 0 } \; . 
\end{equation}
Small deviation asymptotics for Gaussian processes are often reduced  to quantitative estimates on the metric entropy of the unit ball of the Cameron–Martin space when embedded in the given norm. Our study uses a more self contained approach. 
 
The following was shown in \cite{helfter2025scales} as a direct application of Theorem \ref{thm:b-scales} stated below and asymptotics of probabilities of randomly centered small balls by Dereich-Lifshits \cite{dereichproba05}. 
    
    \begin{proposition}
    Let us endow $ \Omega$ with the $ L^p$-norm for $ p \in (1, \infty]$. It holds: 
    \begin{equation*}
     \ord_{ \loc} W ( \bar{\omega} ) = 2 \quad \text{a. s.} 
    \end{equation*}
    Consequently, $ \ord_H W =  2$.
    \end{proposition}
    
  So, any Borel set $A \subset \Omega $ positively supporting $W$ (i.e., $W(A) > 0 $),  has its Hausdorff scale $\ord_H A \ge   2 $ and the equality holds for some Borel set $A$. Therefore, the value of the local order $ \ord_{\loc} W ( \bar{\omega})  $ is equal to $2$ for typical path $ \bar{\omega}$. Here, by ''typical'' we mean ''almost sure'' with respect to the Wiener measure. But there are many ''non typical'' paths.
   The estimation of small ball probabilities for general Gaussian processes is a classical topic, see for instance the survey \cite{LiShao2001}. 
   In most studies, asymptotic probabilities of small balls are given only when they are centered on typical paths, notably  in \cite{dereichproba05} where a wide variety of processes and distance are considered. 
   These results for typical paths are in no way an obstruction for multifractal structure. Non-typical paths are abundant and reveal a realm of multifractal structure of measures on infinite dimension space.  
    We introduce the following space of functions on which we provide the exact multifractal spectrum of the Wiener measure. 
    
    \begin{definition}[Critical Hölder functions]
    For $q \in ( 0  , 1) $, let $C_q^{\crit} $  be the set of \emph{critical $q$-Hölder functions} that is the set of $\bar{\omega} \in \Omega$ such that:
\begin{enumerate}
    \item  the function $ \bar{\omega} $ is $\alpha$-Hölder  for every $\alpha < q $, 
    \item the Hölder exponent of $ \bar{\omega}$ is at most $q$ in the sense that the following limit exists and equals $q$:
    \begin{equation} \label{Def:Orey}
    \lim_{ n \to + \infty } \frac{ - \log \left\{  \sum_{i =1}^n \left[ \bar{\omega} \left(\frac{i}{n}  \right) - \bar{\omega} \left( \frac{i-1}{n} \right) \right]^2 \right\} }{2 \log n } = q \; . 
    \end{equation}
\end{enumerate}
We then introduce the set of \emph{critical Hölder functions} as 
\begin{equation*}
    \Omega^\crit := \bigcup_{ q \in (0,1) }  C_q^{\crit} \;. 
\end{equation*}
\end{definition}
The exponent $q$ given by \eqref{Def:Orey} is introduced in this article and named as  \emph{Orey exponent} of $\bar\omega$.  Note also that $\Omega^\crit $ and all $ C_q^\crit $ are $\sigma$-compact Borel subset of $ \Omega$. We will show in Proposition \ref{prop:ord-loc-fBm} that for 
a typical trajectory of a fractional Brownian motion the Orey exponent coincides with its Hurst parameter. Thus typical trajectories of fractional Brownian motion lie in $ \Omega^\crit $. 
Let us denote: 
  $$ E ( \alpha ) := \lbrace \bar{\omega} \in \Omega  : \ord_{\loc} W ( \bar{\omega} ) = \alpha \rbrace  \, , \qquad f (\alpha) := \ord_H  E ( \alpha)  \;  $$
  and 
    $$ E^\crit ( \alpha ) :=  E ( \alpha) \cap \Omega^\crit  \, , \qquad f_{ \mathrm{crit}} (\alpha) := \ord_H  E^\crit ( \alpha)  \; .   $$
    We will use the convention 
    $ \ord_H \emptyset   = - \infty $. Here is our result on the multifractal spectrum of the Wiener measure $W$.
   
\begin{theo} \label{thm:mf-spec}
    The multifractal spectrum of $W$ on  the space $ \Omega^\crit$  endowed with the uniform norm is fully described by 
    \[ f_\crit ( \xi)  =
    \begin{dcases}
    1 + \frac{\xi}{2}  &  \text{ if } \xi \ge 2  \\
    - \infty   &  \text{ if }  0 < \xi < 2 \;. 
    \end{dcases}
    \]
   Consequently, on the whole space $ \Omega$ it holds
    \[  f ( \xi) \ge   1 + \frac{\xi}{2}
    \qquad \text{ if } \xi \ge 2 \;. \]
    \end{theo} 
    Actually, we will show in Proposition \ref{prop:univ} that $ E (\alpha) = \underline{E}_W^\ord ( \alpha) = \overline{E}_W^\ord( \alpha)  =  \emptyset $ for $ \alpha  < 2 $. A natural question arises whether equalities hold for $ \Omega $ instead of $ \Omega^\crit $ for $ \alpha \ge 2$.  We have no answer for the moment, however in view of Corollary \ref{coro:mass-distrib-pple-intro}, we will simplify this question in Section \ref{lieu:quest-omega-crit}. 
     Note that the spectrum on $ \Omega^\crit$ is not strictly convex, unbounded, with non compact support.
    
    The above Theorem \ref{thm:mf-spec} is implied by our main result  stated below. That is to say,
    we are able to  compute the local order at $\bar\omega  \in C_\alpha^\crit $ of the Wiener measure. 
    \begin{theo}[Local order of Wiener measure] \label{thm:main}
    Let $ \alpha \in (0,1/2)$. Under the $ C^0$-norm, the local order of the Wiener measure $W$ satisfies 
    \begin{equation*}
    \forall \bar{\omega} \in C_\alpha^\crit, \ \ \ 
    \ord_{\loc}  W ( \bar{\omega } ) = 2 ( \alpha^{-1} - 1) \; , 
    \end{equation*}
    while
    \begin{equation*}
    \ord_H C_\alpha^\crit = \alpha^{-1} \; . 
    \end{equation*}
    Moreover the minimum of $ \underline{\ord}_{\loc}  W(\cdot) $ on the space of continuous function $\Omega$ is equal to $2$. 
    \end{theo}
    The proof of Theorem \ref{thm:main}  is decomposed into several intermediary results, most importantly Theorem \ref{thm:ubound}, Theorem \ref{thm:lbound} and Proposition \ref{prop:reg-fBm}\, . 
 The final steps of the proof of Theorem \ref{thm:main}  will be given  in Section \ref{lieu:proofB}, together with the proof of Theorem \ref{thm:fBm-loc-ord} stated just below.
   Theorem \ref{thm:main} suggests us to  use  typical trajectories of fractional Brownian motions to estimate local order of the Wiener measure.
   
     \medskip
     Let us recall that a \emph{fractional Brownian motion} (fBm for short) \( B^H =(B^H_t)_{ t \in [ 0,1]}  \), with \emph{Hurst parameter} \( H \in (0,1) \), is a continuous-time Gaussian process on the segment $ [ 0, 1]$  with initial value \( B^H_0 = 0 \), mean value \( \mathbb{E}[B^H_t] = 0 \) and \emph{covariance} function
    \begin{equation}
     \mathbb{E}[B^H_t B^H_s] = \frac{1}{2} \left( t^{2H} + s^{2H} - |t - s|^{2H} \right) \;. 
    \end{equation}
 See \cite{samoradnitsky2017stable,Kahane1985}  for instance for information about fractional Brownian motions.
   We will prove the following theorem which links Hurst parameter of fBm to the local order of standard Wiener measure at a typical trajectory of the fBm. 
    \begin{theo} \label{thm:fBm-loc-ord}
        Let $ B^H$ be a fractional Brownian motion with Hurst parameter   $ H \in ( 0 ,1)$  and let $W^H$ be the probability law of $B^H$. Then  $W^H$-almost surely the local order  of the standard Wiener measure $W$ (for the uniform norm) at $ B^H$ is given by
        \begin{equation*}
            \ord_{ \loc  } W ( B^H) =
            \begin{dcases}
                2 ( H^{-1} - 1) & \qquad  \text{ if } \, 0 < H < 1/2 \\
                2 & \qquad \text{ if } \, 1/2  \le  H <  1  \, .
            \end{dcases}
        \end{equation*}
    \end{theo} 

\bigskip
    {\em Acknowledgements.} 
    The first author was partially supported by the NSFC (No. 12231013 and No. 12571205) and the
    second author is partially supported by the ERC grant SPERIG $\#885707$.
      The  present work was partially done when the second author visited Wuhan university. 
      The second author would like to thank the
    hospitality of Wuhan university.
    The authors' thanks go to Julien Barral, Arnaud Durand, Shilei Fan and Lingmin Liao for their thoughtful discussions.
    
    \section{A multifractal formalism in infinite dimensional spaces} \label{sec:formalism}
This section lays the foundation for a generalized multifractal formalism in infinite-dimensional spaces,  by introducing rigorously the relevant notions and some basic already known results.
All the notions present in the introduction will be recalled and detailed. 

    \subsection{Scalings}
    We want to study measures defined on a space of infinite dimension.  But the  dimensional multifractal spectrum $f(\alpha)$ is usually infinite. We are going to adopt  
    the setting of scales introduced in \cite{helfter2025scales}, instead of dimensions.
    We consider any metric space $ (X,d)$ and a locally finite Borel measure $\mu$ on it. To present the setting of multifractal analysis in the more general framework , we need to recall some definitions related to scales from \cite{helfter2025scales}. 
    \begin{definition}[Scaling \cite{helfter2025scales}] \label{def:scaling}
        A family $ \Scl = ( \scl_\alpha)_{ \alpha > 0 } $ of Hausdorff functions is  a \emph{scaling} if for every $ \alpha > \beta > 0 $ and every $ \lambda > 1 $ close enough to $ 1 $ (say $1<\lambda <\lambda_0$), it holds: 
        \[ \scl_\alpha ( \varepsilon)  = o ( \scl_\beta ( \varepsilon^\lambda)) \qand \scl_\alpha ( \varepsilon)  = o ( \scl_\beta ( \varepsilon) ^\lambda )  \; .   \]
    \end{definition}
    \begin{remark}
        Note that the supremum value $\lambda_0$ of $ \lambda $ depends on $ \alpha $ and $ \beta $ in the previous definition. 
    \end{remark}
    A particular category of examples of scalings are given by couple of integers  $ p, q \ge 1$ as: 
    \begin{equation} \label{ex:scaling}
    \Scl^{p,q}_\alpha : \varepsilon > 0 \mapsto \frac{1}{\exp ^{\circ p  } ( \alpha  \cdot  \log_+^{ \circ q } ( \varepsilon^{-1 } ) ) }  , 
    \end{equation} 
    for $ \alpha > 0$, where $ \log_+ :=   1_{ (1, + \infty)  } \cdot \log $ is the positive part of the logarithm ($ 1_A $ denoting the indicator function of a set $A$) and $ \circ q$ denotes the $q$-iteration of a map. 
    Note that with $ p = 1 $ and $ q=1$ we retrieve the family classically used to define dimensional invariants.  We will  denote this scaling by $$ \dim :=  \left(  \varepsilon \mapsto \varepsilon^\alpha \right)_{ \alpha>0}.$$ The case of $ p = 2 $ and $  q = 2$ was used to describe some compact spaces of holomorphic maps; see \cite{tikhomirov1993varepsilon,helfter2025scales}. 
    The case of $ p = 2$ and $q=1$  induces what is called the {\em order}. We denote it by $$ \ord  := ( \varepsilon \mapsto \exp( - \varepsilon^{-\alpha} ) )_{ \alpha > 0 }. $$
    The order was used to describe several natural infinite-dimensional objects such as spaces of differentiable maps, see  \cite{tikhomirov1993varepsilon,mcclure1994fractal,helfter2025scales}; ergodic decomposition of measurable maps on smooth manifolds, see \cite{berger2022analytic, berger2017emergence,berger2021emergence,berger2020complexities,helfter2025scales,delaporte2025analytic,berger2025wild} and most importantly for our study here, the geometry of Gaussian processes such as standard Brownian motion \cite{shao1993note,monrad1995small,werner1998existence,dereichproba05}. \medskip
        
    Given a scaling we can define Hausdorff, packing, local scales and then generalize  their corresponding dimensional counterparts. The constructions of Hausdorff and packing measures and their basic properties are yet well known and can be found for instance in \cite{falconer1997techniques,falconer2004fractal,rogers1998hausdorff,tricot1982two}.  
    \subsection{Hausdorff, packing and local scales} \label{lieu:haus-pack}
    For the sake of completeness,
    let us first recall the standard definition of Hausdorff measure. Let $(X,d)$ be a metric space. Let 
    $h$ be a Hausdorff function.

Denote by  $\mathcal{C}_{\rm all}$  the family of all subsets of $X$. 
Given a  sub-family $\mathcal C \subset\mathcal{C}_{\rm all} $ and  a positive number $\delta>0$, for any subset $E\subset X$, define
\[
\mathcal{H}^h_{\mathcal C,\delta}(E)
:=
\inf\Bigl\{\sum_{i=1}^\infty h(\diam A_i)\;:\;
E\subset \bigcup_{i=1}^\infty A_i,\ A_i\in\mathcal C,\ \diam A_i<\delta
\Bigr\},
\]
and then define
\[
\mathcal{H}^h_{\mathcal C}(E)
:=
\lim_{\delta\downarrow 0}\mathcal{H}^h_{\mathcal C,\delta}(E)
=
\sup_{\delta>0}\mathcal{H}^h_{\mathcal C,\delta}(E).
\]
Both $\mathcal{H}^h_{\mathcal C,\delta}$ and $\mathcal{H}^h_{\mathcal C}$ are outer measures; and $\mathcal{H}^h_{\mathcal C}$ restricted on Borel sets is a measure.

If $\mathcal{C}_1, \mathcal{C}_2$ are two families such that $\mathcal{C}_1 \subset \mathcal{C}_2$, it is obvious that 
$$
\forall E \subset X, \quad \mathcal{H}^h_{{\mathcal C}_2}(E) \le \mathcal{H}^h_{{\mathcal C}_1}(E).
$$

We consider the four most common choices  for the sub-family $\mathcal{C}$:
\[
\mathcal C_{\rm all} ,  \qquad
\mathcal C_{\rm open} :=\{G\subset X:\ G\ \text{open}\},
\]
\[
\mathcal C_{\rm closed}: =\{F\subset X:\ F\ \text{closed}\},\qquad
\mathcal C_{\rm ball}: =\{B(x,r): \text{ball}\}
\]
where $B(x,r)$ denotes the open ball centered at $x$ of radius $r>0$. 

Let us  compare  the four associated Hausdorff measures. 
For any $\delta >0$ and $  E\subset X $ the following relationships hold
 \begin{equation} \label{eq:comp-haus-meas}
  \mathcal{H}^h_{\mathcal C_{\rm all}, \delta}(E)=\mathcal{H}^h_{\mathcal C_{\rm closed}, \delta}(E)=\mathcal{H}^h_{\mathcal C_{\rm open}, \delta}(E), \qquad  
   \mathcal{H}^h_{\mathcal C_{\rm all }  , 2 \delta} (E)  \le
\mathcal{H}^h_{\mathcal C_{\rm ball}, 2 \delta}(E)\ \le\ 
\mathcal{H}^{h_2}_{\mathcal C_{\rm all}, \delta}(E) 
\end{equation} 
with $ h_2(t):=h(2t)$. 

We will adopt the most commonly used definition given by all set (or equivalently by open covers) that we will more simply denote as
\begin{equation}
\mathcal{H}_\delta^h := \mathcal{H}^h_{{C_{\rm all}}, \delta}=\mathcal{H}^h_{C_{{\rm open}}, \delta}  \qand \mathcal{H}^h := \mathcal{H}^h_{C_{\rm all}} =\mathcal{H}^h_{C_{\rm open}} \; .  
\end{equation} 
As open covers will be more convenient, by convention, by a $ \delta$-cover of a set $A$ we will mean an open $\delta$-cover of $A$.

 All the above definitions are yet classical and well studied for instance in \cite{rogers1998hausdorff,mattila1997measure,tricot1982two,falconer1997techniques}. This general construction, together with a given  scaling,  allows us  to introduce the notion of Hausdorff scale.
    \begin{definition}[ Hausdorff scale \cite{helfter2025scales}] \label{H scale}
    Given a scaling $ \Scl $. 
     The \emph{Hausdorff scale} of a metric space  $(X,d)$ is defined by: 
    \[\Scl_H  X  := \sup \left\{ \alpha > 0  : \mathcal{H}_{ \mathcal{C}_{ \rm ball }}^{\scl_\alpha}(X) = + \infty  \right\} = \inf \left\{ \alpha > 0  :\mathcal{H}_{ \mathcal{C}_{ \rm ball }}^{\scl_\alpha}(X) = 0   \right\} \; .\]
    \end{definition}
    It is shown in \cite{helfter2025scales} that the second equality  above holds true. Note that thanks to 
 the requirement in the Definition \ref{def:scaling} of scaling and to  the relationship \eqref{eq:comp-haus-meas} we can use
  any of the four subfamilies  $ \mathcal{C}_{\rm ball} $, $  \mathcal{C}_{\rm all}$, $ \mathcal{C}_{\rm open}$ or  $  \mathcal{C}_{\rm closed}$ to define Hausdorff measures and then to get the  same  Hausdorff scale. Indeed replacing $ \scl_\alpha(t) $ by $ \scl_\alpha (2t) $ is harmless in the above definition of Hausdorff scale.

\medskip
We now recall the definition of \emph{packing scale}, following Tricot \cite{tricot1982two}, who defined the packing dimension.
Given $ \delta >0$, a \emph{$\delta$-pack} of a metric space $ (X,d)$  is a countable collection of disjoint balls  of $X$ with radii at most $ \delta$. Let $h$ be a Hausdorff function. As  in the case of Hausdorff outer measure, for $ \delta> 0$, put:   
\[ \mathcal{P}_\delta ^{h}(X) := \sup \left\{ \sum_{i \in I } h (  \radius  B_i   ) : (B_i )_ { i \in I  }\ \text{ is an $ \delta$-pack of $X$}\right\} . \]
    
Since $ \mathcal{P}_\delta^h(X)$ is non-increasing when $\delta$ decreases to $0$,  the following quantity is well defined: 
\[ \mathcal{P}_0^h (X) := \lim_{ \varepsilon  \rightarrow 0} \mathcal{P}_\varepsilon^h(X).\]
As for Hausdorff measures, other definitions of packing measures exist with subtle differences in their behaviours. We adopt the above definition, which  is sufficient for our work. 

\smallskip
 This  allows us to introduce packing measures and packing scales  as follows.
\begin{definition}[Packing  measure \cite{tricot1982two}]For every subset $ E $ of $X$ endowed with the same metric $d$, the \emph{packing  $h$-measure} of $E$ is defined by: 
\[ \mathcal{P}^h (E) = \inf  \left\{ \sum_{n \ge  1 } \mathcal{P}_0^{h} ( E_n ) : E =  \bigcup_{n \ge 1 } E_n  \right\} \; . \]\end{definition}

\begin{definition}[Packing scale \cite{helfter2025scales}] \label{def packing}
The packing scale of a metric space $ (X,d)$ is given by: 
\[ \Scl_P X := \sup \left\{ \alpha > 0  : \mathcal{P}^{\scl_\alpha} (X) = +  \infty \right\}  =\inf  \left\{ \alpha > 0 : \mathcal{P}^{\scl_\alpha} (X) = 0 \right\}   \; . \]
\end{definition} 
    
It is also shown in \cite{helfter2025scales} that the $\sup$ and the $\inf$ here are really equal. The following relation holds. 
\begin{theorem}[Theorem A in\cite{helfter2025scales}]
Let $( X, d) $ be a metric space, then for any scaling $ \Scl $
\begin{equation} \label{ineq:scl}
\Scl_H X \le \Scl_P X \;. 
\end{equation}
\end{theorem} 
Note that it is shown in \cite{hel25} that in any infinite dimensional Banach space and for any scaling, there exists a compact set with arbitrary values of Hausdorff and packing scales as long as \eqref{ineq:scl} holds. 

Now given a Borel measure we recall the definition of local scales. 
\begin{definition}[Local scales \cite{helfter2025scales}] 
Given a scaling $ \Scl $, the upper and lower local scales of a Borel measure $ \mu$ on a metric space $ (X,d) $ are given at a point $x \in X$ by: 
\[ \underline{\Scl}_{ \loc} \mu (x )  = \sup \left\{ \alpha > 0 :   \frac{\mu ( B(x, \varepsilon)) }{\scl_\alpha( \varepsilon) }  \xrightarrow[\varepsilon \to 0]{}  0  \right\} \]
and \[ \overline{\Scl}_{ \loc} \mu (x )  = \inf \left\{ \alpha > 0 :   \frac{\mu ( B(x, \varepsilon)) }{\scl_\alpha( \varepsilon) }  \xrightarrow[\varepsilon \to 0]{}  + \infty   \right\}  \; . \] 
Moreover when they coincide, we denote $ \Scl_{\loc}  ( x) :=  \underline{\Scl}_{ \loc} \mu (x ) = \overline{\Scl}_{ \loc} \mu (x ) $ their common value. 
\end{definition}
        
In particular, if $ \Scl_{ \loc} \mu $ is constant everywhere on $X$, it coincides with both the Hausdorff and packing scale of $(X,d)$.

\begin{remark} \label{rem:supp}
 Observe that if  a point $x \in X$ is outside of the support of $\mu$, meaning $ \mu (B(x,r)) = 0 $ for  small $r$, then  $ \Scl_\loc \mu (x) = + \infty $. This will be used in the formulation of the mass distribution principle. 
\end{remark}

\begin{remark} \label{rem:ex}
Note that when the scaling is of the form $ \Scl^{p,q}$ as in \eqref{ex:scaling}, we have the following expression of local scales:
\begin{equation}
\begin{split}
&  \underline{\Scl}_{ \loc} \mu (x )  = \liminf_{ \varepsilon \to 0 }  \frac{ \log^{\circ p}  \left( \left[  \mu ( B(x, \varepsilon))   \right]^{-1} \right)}{ \log^{\circ q} ( \varepsilon^{-1} ) }   \\
\qand & \overline{\Scl}_{ \loc} \mu (x )  = \limsup_{ \varepsilon \to 0 }  \frac{ \log^{\circ p}  \left( \left[  \mu ( B(x, \varepsilon))   \right]^{-1} \right)}{ \log^{\circ q} ( \varepsilon^{-1} ) }  \; .
\end{split}
\end{equation}
In particular, for the order:
\begin{equation} \label{eq:ord}
\begin{split}
& \underline{\ord}_{ \loc} \mu (x )  = \liminf_{ \varepsilon \to 0 }  \frac{ \log \left( - \log    \mu ( B(x, \varepsilon))   \right)}{  - \log( \varepsilon ) }   \\ \qand & \overline{\ord}_{ \loc} \mu (x )  = \limsup_{ \varepsilon \to 0 }  \frac{ \log \left( - \log    \mu ( B(x, \varepsilon))    \right)}{  - \log( \varepsilon ) } \; . \end{split}
\end{equation}
\end{remark}
While from the probabilistic point of view local scales are of independent interest, in the fractal setting they can serve as a tool to estimates scales of spaces. 
Recall that for  a given scaling $\Scl$, we have defined in \cite{helfter2025scales} the lower and upper Hausdorff scales and the lower and upper  packing scales of a Borel probability measure $\mu$ on a metric space $(X, d)$ respectively by 
\begin{equation}\label{Def:sclH}
\Scl_H \mu  := \inf \left\{ \Scl_H E  :  \mu (E) > 0  \right\} , \quad  \Scl^*_H \mu  := \inf \left\{ \Scl_H E  :  \mu (E) =1  \right\}  ,
\end{equation}
\begin{equation}\label{Def:sclP}  
\Scl_P \mu  := \inf \left\{ \Scl_P E  :  \mu (E) > 0  \right\} , \quad  \Scl^*_P \mu  := \inf \left\{ \Scl_P E  :  \mu (E) =1  \right\} \;  .  
\end{equation}
These scales describe the sizes of sets supporting the measure. 
The following theorem shows that they can be completely characterized by local scales. 
  
\begin{theorem}[Theorem B in 
\cite{helfter2025scales}\label{thm:b-scales}]
Let $ \mu $ be a Borel probability measure on a metric space $(X,d) $. Then  for any scaling $ \Scl$ the following  hold: 
        \[ \Scl_H \mu  
        = \infess  \underline{\Scl}_{ \loc } \mu   , \quad  \Scl^*_H \mu 
        = \supess  \underline{\Scl}_{ \loc } \mu  ,\]
           \[ \Scl_P \mu
           = \infess  \overline{\Scl}_{ \loc } \mu   , \quad  \Scl^*_P \mu 
           = \supess  \overline{\Scl}_{ \loc } \mu   \;  ,  \]
where  $ \supess $ and $ \infess$ denote the essential suprema and infima with respect to  $\mu$.
\end{theorem}
 Such results for Hausdorff dimensions of a measure (i.e. $\Scl_H \mu$ and $\Scl_H^* \mu$  when $\Scl=\dim$), called respectively lower and upper Hausdorff dimensions of $\mu$, were obtained by Fan \cite{fan1994dimensions}. The counterpart of packing dimensions of a measure was later independently considered by  Heurteaux \cite{heurteaux1998estimations} and Tamashiro \cite{tamashiro1995dimensions}.
 
\subsection{Level sets and spectra} 
    Now let us define the multifractal spectrum of a measure relative to a given scaling.
    \begin{definition}[Level sets and multifractal spectrum relative at scale $ \Scl$] 
    The \emph{level sets} of a measure \( \mu \) at scale $ \Scl$  are given for $ \alpha > 0 $ by
    \begin{equation*}
    \underline{E}_\mu^\Scl ( \alpha) := \left\{ x \in X :   \underline{\Scl}_{ \loc} \mu ( x) = \alpha \right\}  \ , \quad 
    \overline{E}_\mu^\Scl  ( \alpha) := \left\{ x \in X :   \overline{\Scl}_{ \loc} \mu ( x) = \alpha \right\} \;  
    \end{equation*}
    and 
    \begin{equation*} 
    E_\mu^\Scl  ( \alpha) := \underline{E}_\mu^\Scl (\alpha) \cap \overline{E}_\mu^\Scl (\alpha)
    =\left\{ x \in X :   \Scl_{\loc} \mu ( x) = \alpha \right\} \; . 
    \end{equation*}
    Then, the (fine) \emph{multifractal spectra} of the measure \( \mu \) at scale $ \Scl$ are defined by the functions
    \[ 
     \underline{f}_\mu^\Scl ( \alpha)  := \Scl_H \underline{E} ( \alpha) 
    \, , \quad 
    \overline{f}_\mu^\Scl ( \alpha)  :=  \overline{E} ( \alpha)   \qand 
   f_\mu^\Scl ( \alpha) := \Scl_H  E ( \alpha)    \; .  \] 
    \end{definition} 
     In the definition of the multifractal spectrum in the dimensional case, it is customary to set $f_\mu^\dim(\alpha)=-\infty$ whenever the level set is empty; and similarly for  $\underline{f}_\mu^\dim$ and $ \overline{f}_\mu^\dim $. 
    This convention marches the multifractal formalism provided  by the Frenchel--Legendre transform (see \eqref{eq:legendre}) and agrees with the expected concavity of the spectrum.  Moreover, it distinguishes genuinely empty level sets from nonempty sets of zero Hausdorff dimension. Consequently we will adopt the same convention for general scales:
    \begin{equation}
        \Scl_H \emptyset := - \infty  \qand \Scl_P \emptyset := - \infty \; .
    \end{equation}

Level sets and multifractal spectra can also be defined using other combinations of local, Hausdorff or packing scales. As an immediate consequence  of Theorem \ref{thm:b-scales}, we obtain the following first step towards multifractal formalism. 
\begin{corollary}[Weak mass distribution principle] \label{cor:spec}
    Let $ \mu$ and $\nu$ be two Borel measures on  a separable metric space $ (X,d) $ and let $ \Scl $  be a scaling. 
    \begin{itemize}
    \item[{\rm (1)}] If $\underline{\Scl}_{\loc}\mu(x) = \alpha$ $\nu$-a.e., then
    $$
      \Scl_H \{x \in X: \underline{\Scl}_{\loc}\mu(x) =\alpha\} \ge \Scl_H^*\nu . 
    $$
    \item[{\rm (2)}] If $\overline{\Scl}_{\loc}\mu(x) = \alpha$ $\nu$-a.e., then
    $$
      \Scl_P \{x \in X: \overline{\Scl}_{\loc}\mu(x) =\alpha\} \ge \Scl_P^*\nu.
    $$
    \end{itemize}
\end{corollary} 
Furthermore, under suitable assumptions, either on the topology of the ambient space or on the form of the scaling functions, we can obtain a reverse inequality yielding a mass distribution principle for the Hausdorff scale of level sets. This will be obtained in the next section in Corollary \ref{coro:mass-distrib-pple} as a direct application of the generalized Frostman Lemma (Theorem \ref{thm:Frostman-gauge}). Note that, for packing measures, no Frostman type result holds even in the dimensional case.

\section{Generalized Frostman Lemma and mass distribution principle} \label{sec-frostman}
 We prove here the Frostman Lemma (Theorem \ref{thm:Frostman-gauge}) and the mass distribution principle (Theorem \ref{thm:mass-distribution-scale}).  An application is to provide a principle which measures the sizes of level sets (Corollary \ref{coro:mass-distrib-pple}).
\subsection{Frostman Lemma for general Hausdorff functions on separable metric spaces} \label{lieu:frostman}

In the theory of Hausdorff dimension in finite dimensional spaces, especially in the Euclidean space $\mathbb{R}^d$, Frostman Lemma plays a fundamental role. We investigate here Frostman Lemma in 
infinite dimensional separable metric spaces.
Let us first state the easy part of Frostman Lemma, whose proof is trivial.
\begin{proposition}[Frostman Lemma, easy part]
\label{thm:Frostman-gauge1}
Let $(X,d)$ be a separable metric space and let $h$ be a Hausdorff function. Consider a Borel set $A$. If there exists a Borel probability measure $\mu$ such that 
\begin{equation}\label{eq:muHolder1}
\mu(A)>0,  \qquad   
\mu\bigl(B\bigr)\, \le\, C \cdot h(\diam B) \ \ \ \forall B ({\rm ball \ in}\ X)
\end{equation}
for some constant $C>0$. Then
\begin{equation}\label{eq:Hausd-estimate}
\mathcal H_{\mathcal{C}_{\text{ball}}}^h(A)\ge \mathcal H^h_{\mathcal{C}_{\text{ball}}, \delta}(A) \ge \frac{\mu(A)}{C}>0 
\end{equation}
for all $0<\delta \le \infty$. 
\end{proposition}

\begin{proof}
For any cover of $A$ by balls $\{B_i\}$ with  diameter  $2r_i\le \delta$,  with $0<\delta\le \infty$ being fixed, we have 
$$
\mu(A) \le  \sum_i \mu(B_i) \le C \sum_i \, h(2 r_i).
$$
Taking the infimum over all such covers 
gives $
\mu(A) \le   C \cdot \mathcal H^h_{\mathcal{C}_{\text{ball}},\delta} (A)
$.  Then, \eqref{eq:Hausd-estimate} follows.
\end{proof}
Note that this result is actually reflected in the first line of Theorem \ref{thm:b-scales}. The conclusion \eqref{eq:Hausd-estimate} remains true for any set $A$ if there exists an outer measure satisfying \eqref{eq:muHolder1}, because only the sub-$\sigma$-additivity is used in the proof. 

 \begin{remark}  We can work with Hausdorff measure defined by arbitrary sets (or equivalently by open sets). We can get a result similar to \eqref{eq:Hausd-estimate},  if $B$ in the condition \eqref{eq:muHolder1}
is assumed to be an open set.  But the condition \eqref{eq:muHolder1} with open set $B$ is not practical. \end{remark}

Frostman \cite{frostman1935potentiel} proved that the inverse of Proposition \ref{thm:Frostman-gauge1} holds when
$X=\mathbb{R}^d$ and $h(t)=t^\alpha$. Kaufman  \cite{kaufman1994min} proved that 
 for  a  general compact metric space $X$, if  the $\alpha$-dimensional Hausdorff measure $\mathcal{H}^\alpha(A)>0$, then for any $\varepsilon>0$ there exists a Borel probability measure $\mu$ on $A$ such that $\mu(B(x,r))\le C r^{\alpha-\varepsilon}$. There is a loss of $\varepsilon$ if we compare with Frostman's original result. But, since $\varepsilon$ is arbitrary, it can deduce that the Hausdorff dimension of $A$ is equal to the capacitary dimension of $A$. This is stated as Frostman Theorem. Thus potential theory is available on compact space to deal with Hausdorff dimension.   
We can improve Kaufman's result by showing that for any gauge $h$ and any $\sigma$-compact metric space, the very Frostman Lemma holds. More generally, we prove the following Frostman Lemma, under the condition that the \emph{increasing sets lemma} holds.
Recall that the gauge $h$ satisfies the increasing sets lemma on a metric space $X$ if, by definition (see \cite{rogers1998hausdorff}, p. 97),  for any increasing sequence of sets $(E_n)$ in $X$ and any couple $0<\delta<\eta$ we have
\begin{equation}\label{eq:increasing-set-lemma}
\mathcal{H}^h_\eta\left(\bigcup_{n=1}^\infty E_n \right) \le \sup_n \mathcal{H}^h_\delta(E_n). 
\end{equation}
As proved by Sion-Sjerve \cite{sion1962approximation} (see also \cite{larman1974}),  this condition is satisfied by  $\sigma $-compact spaces. 
Notice that the space $C^0([0,1])$ of continuous functions defined on 
the interval $[0,1]$ is not $\sigma$-compact. But for any $0<\beta<1$, the space $C^\beta([0,1])$    of $\beta$-H\"{o}lder continuous functions (thus their union as well) is a $\sigma$-compact subspace of $C^0([0,1])$. The Wiener measure describing the Brownian motion is supported by such spaces 
$C^\beta([0,1])$ (with $\beta <1/2$),  to which our Frostman Lemma applies.
Also recall that a subset $A \subset X $ is called \emph{analytic}
when it is the continuous image of a Borel  subset $B$ of some Polish space $Y$, i.e.  there exist a  continuous map \( f : B \to X \) such that \( A = f(B) \).  
In particular, it is well known that every Borel set of a separable metric space is analytic. 
 See  \cite{krantz2008geometric,kechris2012classical,rogers1998hausdorff} for Suslin's theory of analytic sets.

\smallskip

 Let us state and prove 
Frostman Lemma. 
\begin{theo}[Frostman Lemma for a gauge $h$ on a metric space]
	\label{thm:Frostman-gauge}
	Let $(X,d)$ be a separable metric space and let $h:[0,\infty)\to[0,\infty)$ be a Hausdorff function
	(continuous, increasing, $h(0)=0$) such that the increasing sets lemma holds on $X$.
	If $A\subset X$ is analytic and $\mathcal H^h(A)>0$, then there exist a compact set $K\subset A$,
	a Borel probability measure $\mu$ supported on $K$, and a constant $C\in(0,\infty)$ such that
	\[
	\forall B\,  ({\rm ball \ in} \ X), \quad 
	\mu\bigl(B\bigr)\ \le\ C\,h( 3\, \diam B ).
	\]
\end{theo}

Before proving Theorem \ref{thm:Frostman-gauge},  let  us first make a few remarks. 

\begin{remark}
 When $ h (t) = t^\alpha$ with $ \alpha > 0$,  the constant $ 3 $ is absorbed in the constant $ C$ and we retrieve the classical Frostman Lemma and improve Kaufman's result. 
\end{remark} 

 \begin{remark} The condition in Thorem \ref{thm:Frostman-gauge} is made on the Hausdorff measures defined by arbitrary sets and it is stronger than the one made on the Hausdorff measures defined by balls.  In general, 
there are some differences  between these two kinds of Hausdorff measures. \end{remark}

\begin{remark}
One may wonder whether the constant $3$ could be dropped. The answer is negative. In  \cite{davies1969problem},  Davies and Rogers gave  a compact metric space  $( X,d) $ and a  Hausdorff function $h$ such that $\mathcal{H}^h ( A) \in \lbrace 0 , +\infty \rbrace $  for every subset $ A \subset X$. Moreover, they   proved that
every Borel finite measure on $X$ must be singular respectively to $ \mathcal{H}^h$,  namely that they are concentrated on sets of $ \mathcal{H}^h$-measure $0$. But any finite positive measure $ \mu $ such that $ \mu ( B) \le C h  ( \diam B )$ for all balls $ B \subset X$ must be absloutely continuous with respect to  $ \mathcal{H}^h$. 
Yet the question of the optimal constant replacing $3$  is  left open. 
\end{remark}

\begin{question}
 Under the condition \eqref{eq:muHolder1}   in Proposition \ref{thm:Frostman-gauge1}, which holds for balls, can we get 
\eqref{eq:Hausd-estimate} for $\mathcal H^h_{\delta} (A)$  which is defined by arbitrary open sets ? This is an open question, so we can not say that
Theorem \ref{thm:Frostman-gauge} is exactly the converse of Proposition \ref{thm:Frostman-gauge1}.
\end{question} 

\begin{proof}[Proof of Theorem \ref{thm:Frostman-gauge}]
	We will use Theorem~48 (p.~97) in \cite{rogers1998hausdorff}, the Hahn--Banach separation theorem,  the Vitali covering lemma,
	and the Riesz--Markov representation theorem.
The proof  is organized in six steps.
	
	\medskip
	\noindent\textbf{Step 1: Reduction to a compact set with positive $\delta$-content.}
	Recall and denote the $\delta$-level outer measures  defined by open covers as 
	\[ 
	\mathcal H^h_\delta(E)
	:=\inf\Bigl\{\sum_i h( \diam U_i) :\ E\subset \bigcup_i U_i,\ \diam U_i \le \delta\Bigr\},
	\qquad
	\mathcal H^h(E)=\sup_{\delta>0}\mathcal H^h_\delta(E).
	\]
	Note that the analyticity of $ A$ and  the $\delta$-Haudorff  measure $\mathcal H^h(A)$ are unchanged by isometric embeddings into the completion of the metric space. In particular the increasing set lemma holds in the completion of the metric space $(X, d)$. Then we can apply 
 Theorem~48 (p.~97) in \cite{rogers1998hausdorff} to get a compact subset $K \subset A$
	such that $\mathcal H^h(K)>0$. Hence there exists some $\delta\in(0,1]$ such that
	\begin{equation*} 
	 \mathcal{H}_{3\delta}^h ( K ) >0.
	\end{equation*}
	Moreover $\mathcal H^h_{3 \delta}(K)<\infty$ since $K$ is compact and any  open cover admits
	a finite subcover.

	\medskip
	\noindent\textbf{Step 2: A sublinear functional $\Phi_\delta $ on $C(K)$.}
	Let $C(K)$ be the Banach space of real-valued continuous functions on $K$, with norm $\|f\|_\infty:=\sup_{x\in K}|f(x)|$. We are going to construct a non-zero positive sublinear functional on $C(K)$.   
	
Given a non negative function $f \in C(K)$. Consider  a  finite $ \delta$-cover   $\mathcal{B}$ of $K$ by open balls and a strictly  positive function $u:\mathcal{B} \to \R_+:=(0, +\infty)$ such that
 \begin{equation}
  f( x) \le \sum_{ B \in \mathcal{B} \  }  u  (B) 1_B(x)   \;. 
 \end{equation}
 where $1_E$ denotes the indicator function of a set $E$. 
Such a couple $( \mathcal{B}, u)$ will be called a \emph{weighted $\delta$-cover of $f$}.  We will denote by  $ \mathcal{W}_\delta (f) $ the set of weighted $\delta$-cover adapted to $f$.
If $f =1_E$ (the continuity is relaxed), in this way we  formally define a weighted $\delta$-cover of $E$.

For a positive function $f \in C(K)$, we first define
\begin{equation}
\Phi_\delta ( f) := \inf_{ ( \mathcal{B} , u) \in \mathcal{W}_\delta ( f ) }\sum_{ B \in \mathcal{B}} h ( 3\, \diam B ) \cdot u ( B)  \;.
\end{equation}
We then extend $ \Phi_\delta$ as a functional on $C(K)$ by
\begin{equation}
\Phi_\delta ( f) := \Phi_\delta ( \vert f \vert ),  \quad \forall f \in C(K) \;. 
\end{equation}
	Let us now prove the following properties of $ \Phi_\delta$.  Only the proof of the positivity of $\Phi_\delta(1)>0$  is a little bit involved. The proofs of other properties are straightforward.
	
	\smallskip
	{2-a.} {\em $\Phi_\delta$ is bounded}. We have
	\begin{equation}\label{eq:Phi-Lip-upper}
	\Phi_\delta(f)\ \le\ M\,\|f\|_\infty \qquad\forall f\in C(K),
	\end{equation}
	where  $M:=\sum_{B\in \mathcal{B}} h(3\, \diam B )$ for a fixed chosen finite  $\delta$-cover $\mathcal{B}$.  Indeed,   if $\|f\|_\infty \neq 0$ and $u: \mathcal{B}\to \mathbb{R}_+^*$ is constantly equal to $\|f\|_\infty$, then 
	$(\mathcal{B}, u)$ is a weighted $\delta$-cover of $f$. The case $ f = 0$ is similar. 
	
	\smallskip
	 2-b. \emph{$\Phi_\delta$ is positive homogeneous.} We have 
	 \begin{equation}\label{eq:Phi-hom}
	 \Phi_\delta(\lambda f)=\lambda\Phi_\delta(f)\qquad(\lambda\ge 0).
	 \end{equation}
The reason is that  $(\mathcal{B}, u)$ is a weighted $\delta$-cover of $f$ if and only if $(\mathcal{B}, \lambda u)$ is a weighted $\delta$-cover of $\lambda f$  for $\lambda >0$.

	\smallskip
	2-c. \emph{$\Phi_\delta$ is non decreasing on positive functions.}
	Let $f ,g\in C(K)$. Assume $  0 \le f   \le  g $, then 
	\begin{equation}
	\Phi_\delta (f) \le \Phi_\delta(g) \;. 
	\end{equation}
	This is because any weighted $\delta$-cover of $g$ is a  weighted $\delta$-cover of $f$. 
	
	2-d. \emph{$\Phi_\delta$ is subadditivity.} We have 
	\begin{equation}\label{eq:Phi-subadd}
	\Phi_\delta(f+g)\le \Phi_\delta(f)+\Phi_\delta(g)\qquad\forall f,g\in C(K).
	\end{equation}
To prove this, first note that we can assume that $f$ and $g$ are positive,  because
	$$
	 \Phi_\delta ( f+g) = \Phi_\delta ( \vert f + g\vert) \le \Phi_\delta ( \vert f \vert + \vert g\vert)$$
	 where the inequlity is a consequence of  the triangular inequality and the monotonicity (Step 2-c).
	
	We now prove \eqref{eq:Phi-subadd} for positive $f$ and $g$.  Fix $\varepsilon>0$ and 
	choose a weighted $\delta$-covers $(\mathcal{B}_f,u_f)$ and $(\mathcal{B}_g,u_g)$  of $f$  and $g$ such that
	\[
	\sum_{ B\in \mathcal{B}_f} h(3\, \diam B ) \cdot  u_f ( B) \le \Phi_\delta(f)+\varepsilon,
	\qand
	\sum_{ B \in \mathcal{B}_g} h(3\, \diam B ) \cdot  u_g ( B)\le \Phi_\delta(g)+\varepsilon.
	\]
	
Then consider the family $ \mathcal{B} = \mathcal{B}_f \cup \mathcal{B}_g$  (counted with multiplicity) and  define the function $ u: \mathcal{B}\to \mathbb{R}_+^*$ by $ u \vert_{\mathcal{B}_f} = u_f$ and $ u \vert_{\mathcal{B}_g} = u_g $. Then $ ( \mathcal{B} , u) $ defines a weighted $\delta$-cover of $f+g$ since $\mathcal{B}$ forms a $\delta$ cover of $ K$ and for all $x \in K$ it holds 
	\[ 
	 f( x) + g(x) \le \sum_{ B \in \mathcal{B}_f}  u_f (B)1_B(x)   + \sum_{ B \in \mathcal{B}_g }  u_g (B) 1_B(x) = \sum_{ B \in \mathcal{B}}   u (B) 1_B (x)  \;. 
 \]
 It follows then 
	\[
	  \Phi_\delta ( f +g) \le \sum_{ B\in \mathcal{B} } h (3 \, \diam B ) \cdot u ( B ) \le \Phi_\delta ( f) + \Phi_\delta ( f) + 2 \varepsilon \;.  
	\]
	Letting $\varepsilon\downarrow 0$ yields the desired inequality.

	\smallskip
	2-e. \emph{$\Phi_\delta$ is Lipschitz continuous.}  For all $f,g\in C(K)$ it holds
	\begin{equation}\label{eq:Phi-Lipschitz}
	|\Phi_\delta(f)-\Phi_\delta(g)|\le M\|f-g\|_\infty.
	\end{equation} Indeed, 
	from \eqref{eq:Phi-subadd} and \eqref{eq:Phi-Lip-upper}, we get
	\[
	\Phi_\delta(f)\le \Phi_\delta(g)+\Phi_\delta(f-g)\le \Phi_\delta(g)+M\|f-g\|_\infty.
	\]
	Then \eqref{eq:Phi-Lipschitz} follows from this and the symmetry between $f$ and $g$.

	\smallskip
 2-f. \emph{$\Phi_\delta (1) > 0 $. } 
	This last property is the non trivial one and its proof is inspired from Federer \cite{federer2014geometric} (cf. Section 2.10.24).   
	First note that in the definition of $ \Phi_\delta(1)$ we can restrict the infimum over weighted $\delta$-covers  $( \mathcal{B},u)$ such that $u$ takes only rational values. 
	
	Take such a weighted $\delta$-cover $(\mathcal{B}, u)$ of $1_K \in C(K)$. As $ \mathcal{B}$ is finite, there exists a positive integer $N\ge 1$ such that the amplified function $ N \cdot u $ on $ \mathcal{B}$ is  integer valued.  
We will construct recursively a sequence of subfamilies $ \mathcal{B}_0, \dots ,  \mathcal{B}_N $ of $\mathcal{B}$ and non negative integer valued functions $ v_0, \dots , v_N $ starting with $ \mathcal{B}_0 := \mathcal{B} $ and $ v_0 = N \cdot u $ such that for every $ 1 \le j \le N$ the following properties hold
	\begin{enumerate}[(i)]
	\item $\mathcal{B}_j \subset \left\{ B\in \mathcal{B} \ : \ v_{j-1} ( B) \ge 1 \right\} $;
	\item $  \mathcal{B}_j $ is a collection of  
	disjoint open balls of diameter at most $ \delta$, which is called a $\delta /2 $-pack;
	\item $ K \subset \bigcup_{ B (x, r) \in \mathcal{B}_j } B   (x ,3r) $;
		\item $ v_j ( B) = v_{j-1} ( B) - 1 \ge 0 $ if $ B \in \mathcal{B}_j$ and $ v_j (B) = v_{j-1} ( B) $  if $B \in \mathcal{B}\setminus \mathcal{B}_j$; 
		\item  $ K \subset \left\{ x :  \sum_{ B \in \mathcal{B} } v_j ( B) 1_B(x)   \ge N -j  \right\}$.
	\end{enumerate}
To construct this sequence we will proceed inductively. The induction will be similar to the initialization. 

\smallskip
\emph{Initialization.} We first initialize the induction by constructing $ \mathcal{B}_1$
 and $ v_1$. First note that $ \mathcal{B}_0:= \left\{ B \in \mathcal{B} \; : \   v_0 (B) \ge 1 \right\}  $ is the cover $\mathcal{B}$ of $K$ as $u$ is positive.  Thus by Vitali covering lemma for finite covers, it contains a $\delta/2$-pack $ \mathcal{B}_1$ such that the family $ \left\{ B(x,3r) : B(x,r) \in \mathcal{B}_1 \right\}$ covers $ K$. Then $ \mathcal{B}_1$ verifies $ (i)-(iii)$. 

 Note that $(iv)$ defines $v_1$, and $(i)$ ensures that $v_1$ is non negative. 
 Since $\mathcal{B}_1$ is a $\delta$-pack, for every $ x \in K$ there exists at most one ball $ B \in \mathcal{B}_1$ containing $x$. It follows that
 \[ 
 \sum_{ B \in \mathcal{B} } v_1  ( B)1_B(x)  = \sum_{ B \in \mathcal{B}_1 }  (u  ( B)  - 1) 1_B(x) +  \sum_{ B \in \mathcal{B} \backslash \mathcal{B}_1 }   u  ( B)1_B(x)   \ge 
	\sum_{ B \in \mathcal{B} } u  ( B)1_B(x)     \ \ - 1  \, ,\] 
	thus $(v)$ is verified. This concludes the initialization. 
	
	\smallskip 
	\emph{Induction.} Let $ 2 \le j \le N $.   Assume that we have constructed $ v_{j-1} $ and $  \mathcal{B}_{j-1}$ verifying $ (i)-(v)$.  Let us build $ v_j$ and $ \mathcal{B}_j$.
For every $x \in K$, by $(iv)$ for $ j -1$, we have 
\[ \sum_{  B \in \mathcal{B} } v_{ j -1} ( B) 1_B(x)  \ge  N - ( j-1)   > 0  \; \] 
since $ j \le N$.  Then, 
as $v_{j-1}$ is integer valued, there exists $ B \in \mathcal{B}$ such that $ x \in B $ and $ v_{j-1}( B) \ge 1$. This proves that 
 $ \lbrace B \in \mathcal{B}  \ : \ v_{j-1} (B) \ge 1 \rbrace $ covers $ K$. 

By Vitali covering lemma for finite covers,  there exists a subfamily $\mathcal{B}_j \subset \lbrace B \in \mathcal{B}  \ : \ v_{j-1} (B) \ge 1 \rbrace $ verifying $(i)-(iii)$. Then $ (iv)$ defines $v_j$ which is indeed integer valued and non negative by $(i)$.

Since $\mathcal{B}_j$ is a $\delta$-pack, for every $ x \in K$ there exists at most one ball $ B \in \mathcal{B}_j$ containing $x$. It follows that
\begin{align*}
    \sum_{ B \in \mathcal{B} } v_j  ( B) 1_B(x) &= \sum_{ B \in \mathcal{B}_j }  (v_{j-1}  ( B)  - 1) 1_B(x)  +  \sum_{ B \in \mathcal{B} \backslash \mathcal{B}_j }   v_{j-1}  ( B) 1_B(x) \\
   &\ge  \sum_{ B \in \mathcal{B} } v_{j-1}  ( B) 1_B(x)   \ \ - 1   
\end{align*}
 which is at least $ N - (j-1) - 1 = N-j $ by the induction hypotheses, thus $ (v)$ is verified. We then conclude the proof by induction. 
 
\smallskip 
We are now ready to prove $\Phi_\delta(1) >0$. First remark that $\diam \, B(x, 3r)\le 3 \delta$ for each $B(x, r)\in \mathcal{B}_j$. 

By $ (iii) $ which holds for every $ 1 \le j \le N$, we obtain 
\begin{equation} \label{eq:pos}
N \cdot  \mathcal{H}^h_{3 \delta} ( K) \le \sum_{j = 1}^N \sum_{ B\in \mathcal{B}_j} h ( 3\,  \diam B )   \;.
\end{equation}
As  every $ v_{j-1} - v_{j}$ equals $1$ on $ \mathcal{B}_j$ and $ 0$ on $ \mathcal{B}  \, \backslash \,   \mathcal{B}_j$, it follows
\begin{align*}
\sum_{j = 1}^N \sum_{ B\in \mathcal{B}_j} h ( 3 \, \diam B ) &= \sum_{j =1}^N \sum_{ B \in \mathcal{B} }  h (3 \, \diam B)  \cdot  \left[  v_{j-1} (B) - v_{j} (B) \right]  \\
        &= \sum_{ B\in \mathcal{B}} h(3 \, \diam B ) \cdot  \sum_{ j=1}^N \left[ v_{j-1}(B)  - v_j (B) \right] \\
       &= \sum_{ B \in \mathcal{B}} h(3 \, \diam B ) \cdot \left[ N u ( B )  -  v_N(B) \right]\\
     &\le  N   \sum_{ B \in \mathcal{B}} h(3 \, \diam B ) \cdot  u ( B)  \;. 
\end{align*}
This, together with  \eqref{eq:pos}  implies
\begin{equation}
 \mathcal{H}_{3\delta}^h ( K) \le \sum_{ B\in \mathcal{B}} h(3 \, \diam B)  u ( B) \;. 
\end{equation}
Taking the infimum over weighted $\delta$-covers of $1_K$ provides 
$ \Phi_\delta (1) \ge \mathcal{H}_{3\delta} ^h(K)$,  which is strictly positive by  Step $1$.

	\medskip
	\noindent\textbf{Step 3: The single-ball bound.}
	If $0\le f\le 1$ and $\operatorname{supp}f\subset K\cap B$, where $B$ is a ball  of diameter  $\delta$, we have the trivial weighted $ \delta$-cover $ \lbrace \lbrace B \rbrace , 1_{\{B\}}  \rbrace $, then
	\begin{equation}\label{eq:single-ball}
	\Phi_\delta(f)\ \le\ h(3\, \diam B).
	\end{equation}

	\medskip
	\noindent\textbf{Step 4: Hahn--Banach separation and Riesz representation.}
Let us consider 
\begin{equation}
\tau(f)  :=   \inf_K { |f| } \cdot \Phi_\delta (1) = \Phi_\delta ( \inf  { |f| } )    \;.
\end{equation}
Obviously $ \tau$ is concave because $$
\inf_K |\alpha f +(1-\alpha) g| \ge \alpha \inf_K |f| + (1-\alpha)\inf_K |g|.$$
 Moreover  $ 0\le \tau (f)  \le \Phi_\delta  (f) $ for $f \ge 0$. 

Then, by the sandwich version of  Hahn--Banach separation theorem (e.g. \cite{Lassonde1998}, p. 139), there exists a  continuous linear functional $ L : C ( K)  \to \R $ such that
\begin{equation}\label{eq:sandwich}
  L \le \Phi_\delta \ {\rm on }\ C(K)  \qand
	 L \ge \tau  \ {\rm on}\ C_+(K) \; 
\end{equation}
where $C_+(K)$ is the convex set of all non-negative functions in $C(K)$. 
Notice that $ \tau (1) = \Phi_\delta (1) > 0$ by Step 2.  We conclude that the linear functional $L$ is  non zero, positive, continuous and controlled by $\Phi_\delta$. 
	
Consequently, by the Riesz--Markov representation theorem, there exists a finite positive  Borel measure $\mu$ on $K$ such that
	\begin{equation}\label{eq:Riesz}
	L(f)=\int_K f\,d\mu\qquad \forall f\in C(K).
	\end{equation}
	
	\medskip
	\noindent\textbf{Step 5: Estimates for balls of radius $0<r< \delta/2$.}
	Fix $x\in X$ and $0<r<\delta/2$.
	Choose $0<\eta<\delta/2-r$ and, by Urysohn's lemma on the normal space $K$, pick $f\in C(K)$ such that
	\[
	0\le f\le 1,\qquad
	f\equiv 1 \text{ on }K\cap B(x,r),\qquad
	\operatorname{supp}(f)\subset K\cap B(x,r+\eta).
	\]
	Then by \eqref{eq:Riesz}, \eqref{eq:sandwich}, and \eqref{eq:single-ball},
	\[
	\mu\bigl(B(x,r)\bigr)
	\le \int_K f\,d\mu
	= L(f)
	\le \Phi_\delta(f)
	\le h(3\cdot 2(r+\eta)).
	\]
	Letting $\eta\downarrow 0$ and using the continuity of $h$ gives
	\begin{equation}\label{eq:small-scale}
	\mu(B(x,r))\le h(6 \,r)\qquad (0<r<\delta/2).
	\end{equation}
	
	\medskip
	\noindent\textbf{Step 6:  Estimates for balls of arbitrary radius $r>0$ and normalization to a probability measure.}
	For $r\ge \delta/2$, monotonicity of $h$ gives $h(6 \,r)\ge h(3 \delta)>0$, hence
	\[
	\mu(B(x,r))\le \mu(K) 
	    \le \frac{\mu(K)}{h(3 \delta)}\,h(6 \,r).
	\]
	Combining with \eqref{eq:small-scale}, we obtain
	\[
	\mu(B(x,r))\le C_0\,h(6 r)\qquad\forall x\in X,\ \forall r>0,
	\]
	where $C_0:=\max\{1,\mu(K)/h(3 \delta)\}<\infty$.
	
	Finally, set $\nu:=\mu/\mu(K)$. Then $\nu$ is a probability measure supported on $K$ satisfying 
	\[
	\nu(B(x,r))\le C\,h(6 r)
	\qquad\forall x\in X,\ \forall r>0.
	\]
	where $C= \frac{C_0}{\mu(K)}$. This finishes  the proof of the theorem.
\end{proof}

Let us finish this section with some comments on existing investigations related to the Frostman Lemma
and on their links to our proof presented in this paper.

The Frostman Lemma was first proved in Euclidean spaces by Frostman in his thesis \cite{frostman1935potentiel} (1935)   and 
it connects capacities with Hausdorff measures. Frostman's original proof uses the net measures built from covers of dyadic cubes and 
the modern proof of the Frostman uses typically the existence of a set of finite and positive measure.    
The techniques based on  comparable net measures were used by Besicovitch \cite{Besicovitch1952}     
 to obtain subsets of finite positive Hausdorff measure and these techniques were
later formalised by C.\,A. Rogers \cite{rogers1998hausdorff}.  Our proof presented in this paper is based  only on the positivity of $\mathcal{H}_\delta^h(K)$. 

Now let us consider a compact metric space. If the Hausdorff function is of finite order, namely $h(3t) = O(h(t))$ as $t\to 0$,   
the Frostman Lemma follows from Theorem~1 and Theorem~2 in \cite{Howroyd1995}  without passing to a finite-positive subset.
Theorem~2 in \cite{Howroyd1995} is the main result of \cite{Howroyd1995}. A “duality” statement between the weighted Hausdorff measure 
and the measures locally controlled by $h$ (mechanism behind Frostman). Theorem~1 stated in \cite{Howroyd1995}, 
which was essentially proved in Section 2.1.24 of Federer's book \cite{federer2014geometric}, asserts that under suitable conditions, 
the Hausdorff dimension coincides with the weighted Hausdorff dimension when the Hausdorff function
is of finite order.
Note that this property  { being of} finite order is not satisfied in general by the Hausdorff functions in our scalings. 

The weighted Hausdorff measures were first introduced in Kelly's thesis \cite{Kelly1972} (also see his papers \cite{Kelly1973,Kelly1974}). 
 Let us recall the definition of weighted Hausdorff measures.
Let $(X,d)$ be a metric space and let $h:[0,\infty)\to[0,\infty)$ be a Hausdorff function.  
A \emph{weighted cover} of a set $E$ of $X$ is a countable family of pairs $\{(B(x_i,r_i),a_i)\}_{i\in\mathbb{N}}$
such that $a_i\ge 0$ and
\[
1_E(x) \le 
\sum_{i} a_i 1_{B(x_i, r_i)}(x). 
\]
For $\delta>0$, the \emph{$\delta$-level weighted Hausdorff content} is defined by 
\[
\mathcal H^{h}_{w,\delta}(E)
:=\inf\left\{\sum_{i=1}^\infty a_i\,h(r_i)\ :\
\{(B(x_i,r_i),a_i)\}_{i}\ \text{is a weighted cover of }E,\ r_i\le \delta\right\}.
\]
The \emph{weighted Hausdorff measure} is then defined by
\[
\mathcal H^{h}_{w}(E)
:=\sup_{\delta>0}\mathcal H^{h}_{w,\delta}(E)
=\lim_{\delta\downarrow 0}\mathcal H^{h}_{w,\delta}(E).
\]

Actually  one can define $\mathcal{H}^{\xi}_w$ for a more general pre-measure $\xi$  as in \cite{Howroyd1995}. It is clear that our sublinear functional $\Phi_\delta(\cdot)$ is inspired by the weighted Hausdorff measure, a functional variant of  $\mathcal H^{h}_{w,\delta}(\cdot)$. The validity of the increasing sets lemma for ultrametric spaces is due to R.~O.~Davies \cite{Davies1970}.

On compact metric space $X$, a weak version of Frostman
Lemma was proved by R. Kaufman \cite{kaufman1994min}, as we have mentioned. Kaufman used the min-max theorem.  
See \cite{BFP2010} for a detailed proof. 

P. Assouad \cite{Assouad1983} studied  metric spaces which can be bi-Lipschitz embedded into a Euclidean space. On such spaces holds the Frostman lemma for the Hausdorff measure  $\mathcal{H}^h$ with  $h(t) = t^\alpha$  and $\alpha>0$.
According to Assouad \cite{Assouad1983}, a metric space $(X,d)$ is said to be \emph{$(C,s)$-homogeneous} if there exists $C\in(0,\infty)$ and $s\ge 0$ such that, for all $0<a<b$ and for all subsets $Y\subset Z\subset X$ with
\(\text{$Y$ $a$-separated  and $\diam(Z)\le b$},\)
one has the packing bound \[\#Y \ \le\ C\Big(\frac{b}{a}\Big)^{s}.\]
The \emph{metric dimension} of $(X,d)$, now called Assoud dimension, is defined by \[{\rm Dim}_{\rm A}(X,d)\ :=\ \inf\Bigl\{\, s\ge 0:\ \exists\,C\in(0,\infty)\ \text{such that $(X,d)$ is $(C,s)$-homogeneous}\Bigr\}.\]
We say that $(X,d)$ has \emph{finite metric dimension} if ${\rm Dim}_{\rm A}(X,d)<\infty$.
Assouad proves an embedding theorem (Proposition 2.6 in \cite{Assouad1983}): if a metric space $(X,d)$ has finite metric dimension and if $0<p<1$, then $(X, d^p)$ admits a bi-Lipschitz embedding into some $\mathbb{R}^n$. 
Thus the Frostman lemma holds for spaces of finite metric dimension. This result due to Assouad is now covered by Howroyd's result because the Hausdorff function $t\mapsto t^\alpha$ is of finite order. We point out that the finite metric dimension condition of $(X, d)$ is equivalent to the doubling property of $(X, d)$, namely, there exists an integer $N\ge 1$ such that for every $x\in X$ and $R>0$, the ball $B(x,R)$ can be covered by at most $N$ balls of radius $R/2$.

\subsection{Mass distribution principle and its application to level sets}
From the Frostman Lemma we have just proved,  it is easy  to obtain the following mass distribution principle at any scale. That is to say that, the Hausdorff scale of an analytic set is equal to the supremum of the Hausdorff scales of the measures concentrated on the set. 
Notice that we can always assume that the Borel measure $\mu$ involved in Theorem \ref{thm:mass-distribution-scale} below is complete, and that analytic sets are $\mu$-measurable, so that $\mu\in \mathcal{M}^1(A)$ makes sense.

\begin{theo}[Mass distribution principle]
\label{thm:mass-distribution-scale}
Let $ 
 \Scl = ( \scl_\alpha)_{ \alpha > 0 } $ be a scaling.
 Let $(X,d)$ be a separable metric  such that each $ \scl_\alpha $ satisfies the increasing sets lemma on $X$. 
Then for any analytic subset  $ A \subseteq X$ the Hausdorff scale verifies the mass distribution principle
\[
\Scl_H A = \sup_{\mu \in \mathcal{M}^1(A)} \Scl_H^* \mu =  \sup_{\mu \in \mathcal{M}^1(A)} \Scl_H \mu = \sup_{\mu \in \mathcal{M}^1(A)}  \inf \underline{\Scl}_\loc \mu,
\]
where $\mathcal{M}^1(A)$ denotes the set of Borel probability measures supported on $A$. 
\end{theo}

\begin{proof} 
First note that by definitions of Hausdorff scales of measures
\begin{equation}
 \sup_{\mu \in \mathcal{M}^1(A)} \Scl_H^* \mu \ge  \sup_{\mu \in \mathcal{M}^1(A)} \Scl_H \mu \ge \sup_{\mu \in \mathcal{M}^1(A)}  \inf_{ x \in A} \underline{\Scl}_\loc \mu  (x) \; . 
\end{equation}
Moreover they are all at most $ \Scl_H A  $ by Proposition \ref{thm:Frostman-gauge1} (or the first line in Theorem \ref{thm:b-scales}). Thus, to conclude, it remains only to show that
\begin{equation}
  \Scl_H A  \le    \sup_{\mu \in \mathcal{M}^1(A)}  \inf_A\underline{\Scl}_\loc \mu  \;.
\end{equation}
Fix $ \alpha  < \Scl_H A$. By the definition of Hausdorff scale, $\mathcal H^{\scl_\alpha}(A)  = +\infty$. Thus, by  Frostman Lemma (Theorem \ref{thm:Frostman-gauge}), 
there exists a probability measure $\mu$ supported on $A$ and a constant $C>0$ such that 
\[
\mu(B(x,r)) \le C \cdot \scl_\alpha(6 r), \quad \forall x\in A, \, r>0.
\]
By Definition \ref{def:scaling} of scaling, for every $ \beta > \alpha$ and $x \in  A $ it follows
\[
\mu(B(x,r))  = o  \left(   \scl_\beta (r) \right) \quad \text{ as } r \to 0 \;. 
\] 
Consequently it holds that 
$
\underline{\Scl}_\loc \mu(x) \ge \alpha
$ for every $x \in A$. 
 In particular
$$
\sup_{\nu \in \mathcal{M}^1(A)}  \inf_{ x \in A} \underline{\Scl}_\loc \nu  (x)   \ge \alpha.
$$
Now we can conclude because $\alpha <\Scl_H A$ is arbitrary.
\end{proof}
Note that we have thus proved in particular the result of Theorem \ref{thm:mass-distribution-scale-intro}. 
\medskip

Then to apply the latter Theorem \ref{thm:mass-distribution-scale} to level sets we shall show that they are analytic. In fact, they are Borel sets. We first prove the following lemma. 
\begin{lemma} \label{lem:borel}
Let $ (X,d) $ be a separable metric space and $ h : \R_+ \to \R_+ $ be a non decreasing function, then 
\[  \underline{G} :=\left\{x\in X:\limsup_{r\to 0}\frac{\mu(B(x,r))}{h(r)}=0\right\} \]
and \[ 
\overline{G} :=\left\{x\in X:\liminf_{r\to 0}\frac{\mu(B(x,r))}{h(r)}= + \infty \right\} 
\]
are Borel sets. 
\end{lemma}
\begin{proof}
For fixed integers \(n,m\) define 
\[
\underline{A}_{n,m}:=\bigcap_{\substack{q \in \Q \\ 0<q<1/m}}
\Big\{x\in X \, :\, \mu(B(x,q))<\tfrac{1}{n}\,h(q)\Big\} \]
and
\[
\overline{A}_{n,m}:=\bigcap_{\substack{q \in \Q \\ 0<q<1/m}}
\Big\{x\in X \, :\, \mu(B(x,q))> n \,h(q)\Big\}.
\]  
 which  are Borel measurable since the maps  $ (x,r) \mapsto \mu (B(x,r)) $  and $h$ are Borel measurable.
Then, as $h$ is non decreasing, we can write 
\[
\underline{G} =\bigcap_{n\ge1}\bigcup_{m\ge1} \underline{A}_{n,m} \qand 
\overline{G} =\bigcap_{n\ge1}\bigcup_{m\ge1} \overline{A}_{n,m}
\] 
which are countable intersections of countable unions of Borel sets.
\end{proof}
Consequently, we obtain measurability of level sets. 
\begin{proposition} \label{prop:lvl-borel}
Let $ \Scl $ be a scaling and $ \mu $ a Borel measure on a separable metric space $ (X,d) $. Then for every $\alpha > 0$ the  level sets
 \( \underline{E}_\mu^\Scl(\alpha) ,  \overline{E}_\mu^\Scl(\alpha)  \)  and  \( E_\mu^\Scl(\alpha)  \) are Borel thus analytic subsets of $ X$.  
\end{proposition}
\begin{proof}
    For $ \beta > 0 $, let 
    $$ \underline{G}( \beta ) := \left\{x\in X:\limsup_{r\to 0}\frac{\mu(B(x,r))}{\scl_\beta(r)}=0\right\} $$ and 
    $$\overline{G}( \beta ) := \left\{x\in X:\liminf_{r\to 0}\frac{\mu(B(x,r))}{\scl_\beta(r)}= + \infty\right\}   \;. $$ 
 By the above Lemma \ref{lem:borel}  these sets are  Borel measurable, then observe that the level sets $ \underline{E}_\mu^\Scl(\alpha) $ and $\overline{E}_\mu^\Scl(\alpha) $  write as:  
    \begin{equation*}
         \underline{E}_\mu^\Scl(\alpha)  = \bigcap_{ \substack{ \beta \in \Q    \\ 0 < \beta < \alpha } } \underline{G}( \beta)  \backslash \bigcup_{ \substack{ \beta \in \Q    \\  \beta > \alpha } } \underline{G}( \beta)  \qand \overline{E}_\mu^\Scl(\alpha)  = \bigcap_{ \substack{ \beta \in \Q    \\  \beta > \alpha } }\overline{G}( \beta)  \backslash \bigcup_{ \substack{ \beta \in \Q    \\ 0 < \beta < \alpha } }  \overline{G}( \beta)  
    \end{equation*}
    are  countable combinations of Borel sets. It follows also that $ E_\mu^\Scl(\alpha)  = \underline{E}_\mu^\Scl(\alpha)  \cap \overline{E}_\mu^\Scl(\alpha) $ is Borel measurable.
\end{proof}

As a direct application of Theorem \ref{thm:mass-distribution-scale} and the above Proposition \ref{prop:lvl-borel} we obtain the following. 
\begin{corollary} \label{coro:mass-distrib-pple}
    Let $ \mu$ be a Borel measure on  a separable metric space $ (X,d) $ and let $ \Scl $  be a scaling such that $ \scl_\alpha $ verifies the increasing sets lemma for every $ \alpha >0$. Then for every $ \alpha > 0$,  the Hausdorff scales of the level  sets $   \underline{E}_\mu^\Scl(\alpha)$ and $  E_\mu^\Scl(\alpha)$ are respectively  equal to
      $$    \underline{f}_\mu^\Scl(\alpha)  =   \sup \left\{ \Scl_H \nu \, : \, \nu \in \mathcal{M}^1 \left( \underline{E}_\mu^\Scl(\alpha) \right)   \right\} =  \sup \left\{ \Scl_H^* \nu \, : \, \nu \in \mathcal{M}^1 \left( \underline{E}_\mu^\Scl(\alpha) \right)   \right\} \;      $$
      and 
            $$    f_\mu^\Scl(\alpha)  =   \sup \left\{ \Scl_H \nu \, : \, \nu \in \mathcal{M}^1 \left( E_\mu^\Scl(\alpha) \right)   \right\} =  \sup \left\{ \Scl_H^* \nu \, : \, \nu \in \mathcal{M}^1 \left( E_\mu^\Scl(\alpha) \right)   \right\} \;    .  $$
\end{corollary}

    \section{Multifractal analysis of the Wiener measure} \label{sec:wiener}
 
In this section, we compute the multifractal spectrum of the Wiener measure. That is to say,  we are going to prove Theorem \ref{thm:mf-spec}. We first show Theorem \ref{thm:fBm-loc-ord} concerning  the local scales of the fBm's  and Theorem \ref{thm:main} concerning 
 the scales of the sets of critical Hölder functions. Theorem \ref{thm:mf-spec}
will easily follows from Theorem \ref{thm:fBm-loc-ord}  and Theorem \ref{thm:main}.

The proof of the lower bound in Theorem \ref{thm:main} uses a fine analysis of increments and the Orey exponent (Proposition \ref{prop:reg-fBm}).
The proof of the upper bound uses the cylinder covering  (Theorem \ref{thm:cyl-ball}).
Combining these yields Theorem \ref{thm:main}.
Our approach combines estimates on small ball probabilities based on  fine geometric analysis of the local behavior of trajectories.
We begin by recalling the asymptotics of the Wiener measure of small balls due to Dereich and Lifshits \cite{dereichproba05}, which describe the exponential decay of the measure of typical small balls in the uniform norm. 
As a first step, 
we show that the lower bound in their result extends to the entire space of continuous trajectories (Proposition \ref{prop:univ}). 
Therefore, the level sets of local order strictly less than \(2\) are empty. Next, we relate the local order of the Wiener measure to the regularity of the underlying path (i.e., the center of the ball) through its Hölder exponent.  An upper bound  of local order of the Wiener measure is obtained in Theorem \ref{thm:lbound}, and  a lower bound is obtained in Theorem \ref{thm:ubound}. The lower bound is described by the notion of Orey exponent of the considered trajectory which is introduced here in \eqref{def:orey}. 
The proofs of both results rely crucially on the fact that Brownian motion has independent stationary increments.  This allows to effectively estimates probabilities of balls by comparing them with the probabilities of  suitably chosen cylinders.  
Combining these bounds shows that, for trajectory $\bar\omega$ belonging to the critical Hölder space \( \Omega^\crit \) then to some \(C_\alpha^\crit \), 
the local order  of the Wiener measure at $\bar\omega$ is well defined and equal to \( 2(\alpha^{-1} - 1) \).
Classical regularity results for fBm, that can be found in \cite{Xiao2009} by Xiao, yield that typical trajectories of fBm with Hurst parameter \(H\) have Orey exponent equal to their maximal Hölder regularity \( H \) (Proposition \ref{prop:reg-fBm}). 
To conclude the proofs of Theorem \ref{thm:main} and Theorem \ref{thm:fBm-loc-ord}, we
compute the Hausdorff orders of different spaces of Hölder functions. In particular, using the mass distribution principle, we prove that the Hausdorff orders of fBms are the Hausdorff orders of the corresponding level sets restricted to \( \Omega^\crit \). This is achieved using the  estimates of measures of balls centered at the origin obtained by Li and Werner \cite{werner1998existence}, together with the general Gaussian process results of Dereich and Lifshits \cite{dereichproba05}, which extend these estimates from the origin to almost all trajectories. These two key results are reformulated in the framework of scales in Proposition \ref{prop:ord-loc-fBm}, leading to the partial classification of non-typical trajectories stated in Theorem \ref{thm:main} and Proposition \ref{prop:reg-fBm}.

Throughout this section, the estimates on small ball probabilities are interpreted geometrically in terms of scales, providing a natural bridge between the small-deviation theory of Gaussian processes and the infinite-dimensional multifractal analysis.

    \subsection{Local order of Wiener measure and Orey exponents of fractional Brownian motions} \label{lieu:statement-main}
    
    Let us recall that $ \Omega := \left\{  \omega \in C^0( [0,1] , \R ) : \omega(0) = 0 \right\}$, which is endowed with the supremum norm and that a function $ \omega \in \Omega$ is  said to be $ \alpha$-Hölder for $ \alpha \in  (0, 1]$ when the $ \alpha$-Hölder modulus of $\omega$ is finite, i.e.  
    \begin{equation}
     \| \omega \|_\alpha := \sup_{ x \neq y}  \frac{ \vert \omega(x)-\omega(y) \vert}{ \ \  \vert x-y \vert^{ \alpha }  } \in [ 0 , + \infty ) \;  . 
    \end{equation} 
    We will be interested in the subspaces $C_\alpha$  of $\alpha$-H\"{o}lder functions: 
    \begin{equation}
        C_\alpha := \left\{ \omega \in \Omega : \| \omega \|_\alpha < + \infty \right\} , \quad \alpha \in (0,1]. 
    \end{equation}
    We will also need the space of 
    \emph{nearly $\alpha$-Hölder} functions:
    \begin{equation}
     C_{\alpha}^-:= \bigcap_{0<\beta<\alpha} C_\beta, \quad \alpha \in (0,1].
    \end{equation}
    It is well known that the trajectories of Brownian motion lie almost surely in the space $ C_{1/2}^-$ and the exponent $1/2$ is optimal. More generally, a fBm with Hurst index $ H \in (0,1) $ has almost surely its trajectories in $ C_H^-$. This can be found for instance in \cite{Kahane1985} (p.264). 
 As we have already pointed out in the introduction,  we will describe the multifractal spectrum  of the Wiener measure $W$ by studying  the local order $\ord_{\loc}W(\bar{\omega})$  at the typical trajectories $\bar{\omega}$ of fBm's (see Theorem \ref{thm:fBm-loc-ord} and Theorem \ref{thm:main}).

    We introduce a notion of maximal Hölder exponent that we call \emph{Orey exponent}. 
    We first define the discrete quadratic variation of a function $ \bar{\omega} \in \Omega$ given for $ n \ge 1$ by 
    \begin{equation}
    Q_n(\bar{\omega}) :=   \sqrt{ \frac{1}{n} \sum_{i=1}^{n} \left[ \bar{\omega}\left( \frac{i}{n} \right) - \bar{\omega}\left( \frac{i-1}{n} \right) \right]^2 }  \;.
    \end{equation}
    Then we define the \emph{lower and upper Orey exponents} respectively by
    \begin{equation} \label{def:orey}
        q^- ( \bar{\omega} ) := \liminf_{ n \to + \infty} \frac{-\log Q_n ( \bar{\omega} ) }{  \log n } \qand   q^+ ( \bar{\omega} ) := \limsup_{ n \to + \infty} \frac{-\log Q_n ( \bar{\omega} ) }{  \log n } \;. 
    \end{equation}
    We say that $ \bar{\omega} $ has \emph{Orey exponent}  $q (\bar{\omega}) $ if  $  q( \bar{\omega} ) :=  q^- ( \bar{\omega} )  =  q^+ ( \bar{\omega} ) $.
    This Orey exponent captures in some $ L^2$-sense the maximal Hölder exponent of the considered trajectory. 
    Indeed, for $\bar{\omega}\in C_\alpha$, we clearly have $Q_n(\bar\omega)=O(n^{-\alpha})$ so that $q^-(\bar{\omega})\ge \alpha$, namely $q^-(\bar\omega)$ upper bounds the H\"{o}lder exponent of $\bar\omega$. 
 Let us state clearly this fact. 
    \begin{fact} \label{fact:mod}
    If $ \bar{\omega}  \in C_\alpha^-$, then
    $q^- ( \bar{\omega} ) \ge  \alpha $.
    \end{fact}
     The introduction of the Orey exponent is motivated by the notion of Orey index which was proposed by Orey in \cite{orey1970gaussian} to characterize Hölder continuity properties of Gaussian processes such as fBm. The following proposition says that almost surely, the trajectories of a fBm of Hurst index $H$ are nearly $H$-H\"{o}lder and has $H$ as its Orey exponent. The proof will use the self-similarity of $B^H$: 
\[
  \forall c>0, \quad (B^H_{ct})_{t\in\mathbb{R}} \stackrel{d}{=} (c^H B^H_t)_{t\in\mathbb{R}}
\]
and the stationary  of increments:
\[
\forall s,t,h\in\mathbb{R}, \quad 
  B^H_{t+h} - B^H_{s+h} \stackrel{d}{=} B^H_t - B^H_s.
\]

Let us also recall the following useful fact.
The fractional Brownian motion $B^H=(B_t^H)_{t\in\mathbb R}$ can be defined by the 
Mandelbrot--Van Ness formula (see \cite{MvN1968}):
\[
B_t^H = \int_{\mathbb R} K_H(t,u)\,dW_u
\]
with
\[
K_H(t,u)= \frac{1}{ \Gamma ( H +1/2 ) }  \Big((t-u)_+^{H-\frac12}-(-u)_+^{H-\frac12}\Big),
\]
where $(W_u)_{u\in\mathbb R}$ is a two-sided standard Brownian motion  and $ x_+ := \max ( 0, x) $. 
The normalizing constant $\Gamma ( H +1/2 )  $, where $ \Gamma $ is the usual Gamma function, is chosen so that
$\E[(B_t^H)^2]=|t|^{2H}$.
Fix $s<t$. Set $$
\Delta_{s,t}B^H := B_t^H-B_s^H.
$$
Then the centered random variable
\[
(\Delta_{s,t}B^H)^2 - \E\big[(\Delta_{s,t}B^H)^2\big]
\]
belongs to the second Wiener chaos 
$\mathcal{H}_2(W)$ of $W$. 
To this centered variable, we will apply the  
hypercontractivity (Nelson) inequality of Wiener chaos. We refer to \cite{NourdinPeccati2012} for the Wiener chaos and its hypercontractivity.

    \begin{proposition}[Orey exponent of fBm] \label{prop:reg-fBm}
    Let $ B^H $ be a fractional Brownian motion with Hurst parameter $ H \in (0,1)$. Then
    \begin{equation*}
    q ( B^H)=H  \qand B^H \in  C_{H}^-  \quad  W^H-\text{a.s} 
    \end{equation*} 
    \end{proposition} 
    \begin{proof}
    The result
    $B^H \in C_H^-$ on the modulus of continuity of $B^H$ is  well known, see \cite{Kahane1985,marcus2006markov}. 
    Actually, we have the following precise result, called
    the sharp uniform modulus of continuity of fBm (L\'evy-type law):
\begin{equation}\label{eq:modulus}
    \limsup_{r\downarrow 0}\ \frac{\sup_{|t-s|\le r}|B^H_t-B^H_s|}
    {r^{H}\sqrt{2\log(1/r)}}=1\qquad\text{a.s.}
    \end{equation}
    (see  Xiao \cite{Xiao2009}).  
    This implies  that $B^H$ is a.s. $\alpha$-H\"{o}lder for every $\alpha <H$. Then, by using Fact \ref{fact:mod} and the proved fact $B^H\in C_H^-$, we get immediately the lower bound $q^-(B^H)\ge H$ a.s. Thus it remains to prove the upper bound $q^+(B^H)\le H$ a.s.

    \smallskip
    Consider the $1/n$-mesh increments
    \[
    \Delta^{(n)}_j\ :=\ B^H\!\left(\tfrac{j+1}{n}\right)-B^H\!\left(\tfrac{j}{n}\right),\qquad j=0,\dots,n-1.
    \]
    By the self-similarity and the stationary of increments of $B^H$,
    \begin{equation}\label{eq:law-eq}
    \left(\Delta^{(n)}_0,\dots,\Delta^{(n)}_{n-1}\right)\ \stackrel{d}{=}\ n^{-H}\,(\xi_0,\dots,\xi_{n-1}),
    \end{equation}
    where  $\xi_j:=B^H_{j+1}-B^H_j$  for $0\le j \le n-1$ is a
     \emph{fractional Gaussian ``noise''}, namely, a centered stationary Gaussian sequence with
    \[
    \mathbb{E}[\xi_0^2]=\mathrm{Var}(B^H_1)=1,\qquad 
    \mathbb{E}[\xi_0\xi_k]=\rho_H(k):=\tfrac12\left(|k+1|^{2H}-2|k|^{2H}+|k-1|^{2H}\right).
    \]
    
Observe that,   if $ H \neq \frac{1}{2}$, then  $\rho_H(k)\sim H(2H-1)\,   \vert  k \vert ^{2H-2}$, and if $ H = 1/2 $ then $ \rho_h( k) = 0 $ for $  k \neq 0 $.  
Hence, for every $ H \in ( 0,1)$  there exists a constant $\kappa_H > 0 $ such that: 
 \begin{equation}\label{eq:xi-k}
          \vert \rho_H ( k ) \vert \le \kappa_H \vert k \vert^{2H -2 }  \text{ for }  k   \neq 0 \;. 
 \end{equation}
 In particular $\rho_H(k)\to 0$
as $k\to\infty$ since $ H < 1$. 
 Let us introduce the normalized quadratic variation
    $$
    V_n := n^{2H}Q_n(B^H)^2
    =\frac{1}{n}\sum_{j=0}^{n-1}\left(n^{H}\Delta^{(n)}_j\right)^2.
    $$
    From \eqref{eq:law-eq}, we have
     \[
    V_n
    \ \stackrel{d}{=}\ \frac{1}{n}\sum_{j=0}^{n-1}\xi_j^2 =:W_n.
    \]
    We are going to show that  $\sum\mathbb{P} \left(V_n  \le 1/2  \right) <\infty$ with $\mathbb{P}=W^H$,   so that we can apply Borel-Cantelli lemma. Since $V_n$ and $W_n$ obey the same probability law, we will work with $W_n$. First it is obvious that $\mathbb{E} ( W_n ) = 1$.  
    
    \smallskip
  Secondly, we provide an estimate of the variance of $W_n$,  using a standard identity for centered Gaussian variables. If $U,V$ are centered Gaussian with $\mathbb E[U^2]=\mathbb E[V^2]=1$ and
$\mathbb E[UV]=\rho$, then
\[
\operatorname{Cov}(U^2,V^2)=2\rho^2.
\]
This is a special case of 
Isserlis’ theorem or  Wick’s formula for Gaussian moments \cite{isserlis1918formula}, see for instance \cite{munthe2025short} for a short and general proof.
Applying  the above identity to the pairs $(\xi_j^2, \xi_k^2)$ and using  the stationarity of $(\xi_k)$, we obtain
\[
\begin{aligned}
\operatorname{Var}(W_n)
&=\frac1{n^2}\sum_{j,k=0}^{n-1}\operatorname{Cov}(\xi_j^2,\xi_k^2) 
 =\frac2{n^2}\sum_{j,k=0}^{n-1} \left[\rho_H(k-j) \right]^2 \\
 &=\frac2{n^2}\Bigl(
  n\,\rho_H(0)^2 + 2\sum_{\ell=1}^{n-1}(n-\ell)\rho_H(\ell)^2\Bigr).
\end{aligned}
\]
Using $\rho_H(0)=1$ and \eqref{eq:xi-k},
we obtain
\begin{equation} \label{eq:var-wn-1}
\operatorname{Var}(W_n)
\le \frac{2}{n} + \frac{4\kappa_H}{n^2}
      \sum_{\ell=1}^{n-1}(n-\ell)\,\ell^{4H-4}.
\end{equation}
 We use the trivial 
 estimation 
\begin{equation} \label{eq:triv}
\sum_{\ell=1}^{n-1}(n-\ell)\,\ell^{4H-4} 
\le  n \sum_{\ell=1}^{n-1} \ell^{4H-4} \;. 
\end{equation}
Then note that 
\begin{equation*}
   \begin{dcases}
   \sum_{ \ell \ge 1} \ell^{4H-4} < + \infty  & \qquad  \text{ if } \,   0 < H < 3/4  \\
  \sum_{ \ell = 1}^{n-1} \ell^{4H-4}  = O (  \log n )    & \qquad  \text{ if } \,    H  = 3/4   \\
  \sum_{ \ell = 1}^{n-1} \ell^{4H-4} =  O (  n^{4H - 3}  )    & \qquad  \text{ if } \,    3/4 < H < 1    \\
  \end{dcases}
  \end{equation*}
   as $ n \to + \infty$. Considering all cases we obtain \begin{equation} \label{eq:sum}
    \sum_{\ell= 1}^{n-1} \ell^{4H - 4} = O \left(   \log n  +  n^{4H-3} \right)  \;. 
\end{equation}   
 Thus,   if we take  $ 0<\delta <\min\{1,   4(1-H)\}  $, the estimations  \eqref{eq:var-wn-1}, \eqref{eq:triv} and \eqref{eq:sum} together imply the existence of a constant $ \kappa_{\delta,H}$ such that 
\begin{equation} \label{eq:var-wn}
    \forall n \ge 1, \qquad \operatorname{Var}(W_n)\le   \kappa_{\delta,H}  \cdot n^{-\delta}.
\end{equation}
We now show that for every $p\ge 2$, there exists $ C_p $ such that  \(\mathbb E|W_n-1|^p \le \, C_p\,n^{- \delta p/2} \). 
From
$$ W_n-1 = \frac{1}{n}\sum_{j=0}^{n-1}  (\xi_j^2-1) \; , 
$$
we see that 
the variable $W_n-1$ belongs to the second Wiener chaos of the underlying
Gaussian space. As an application of the hypercontractivity (Nelson) inequality, it satisfies for every $ p \ge 2$ 
\[
 \left[ \E|W_n-1|^p \right]^{1/p} 
\le C_{p,2}   \,\left[ \Var(W_n) \right]^{1/2} \; 
\]
where $ C_{p,2} > 0 $ is a universal constant, see e.g.\ Hairer~\cite[Theorem 7.1]{hairer2021introduction}.
The latter inequality together with \eqref{eq:var-wn} provides the existence of a constant $C_p > 0  $ (depending on $H$ and $\delta$ as well) such that
\begin{equation} \label{eq:est-v_np}
    \mathbb E|V_n-1|^p =  \mathbb E|W_n-1|^p
\le C_p \,n^{-\delta p/2}.
\end{equation}
Note  then that  $ \{V_n\le\tfrac12\}
\subset\{|V_n-1|\ge\tfrac12\}  $. Markov's inequality and the estimate in \eqref{eq:est-v_np}  provide
\[
\mathbb P\bigl(|V_n-1|\ge\tfrac12\bigr)
\le 2^p\,\mathbb E|V_n-1|^p
\le 2^p C_p \,n^{- \delta p/2}.
\]
Choosing $p$ large enough so that $\delta p/2>1$, we obtain
\[\sum_{n=1}^\infty \mathbb P\bigl(V_n \le\tfrac12\bigr) \le
\sum_{n=1}^\infty \mathbb P\bigl(|V_n-1|\ge\tfrac12\bigr)
<\infty.
\]
Hence, by the Borel--Cantelli lemma,
\(
\mathbb P\bigl(V_n\le  \frac{1}{2} \, \text{i.o.}\bigr)=0.
\)
Equivalently, a.s. $Q_n(B^H) > \frac{1}{\sqrt{2}} n^{-H}$ for large $n$, which implies $q^+(B^H)\le H$ a.s.
    \end{proof}
This leads us to consider the space of continuous functions with maximal H\"{o}lder exponent $\alpha$,  for $ \alpha \in (0,1]$: 
     \begin{equation}
     C_\alpha^\crit := \left\{ \omega \in C_\alpha^- : q(\omega) = \alpha \right\} \; .
     \end{equation}
     Note that Theorem \ref{thm:main} is actually equivalent to 
     \begin{equation}
          C_\alpha^\crit \subset E_W^\ord \left(  \max \left[ 2 , 2 ( \alpha^{-1} - 1 )  \right] \right) \qand \ord_H  C_\alpha^\crit  = \alpha^{-1}
     \end{equation}
     for every $ \alpha \in ( 0,1)$. 
   To prove this we first show that asymptotics of probabilities of  small balls can be obtained by computing the measure of some well chosen cylinders. 
   
    \subsection{From small ball estimates  to cylinders estimates} \label{lieu:ball-to-cyl}
    
    As the Wiener measure has independent increments,  it is quite simple to obtain an estimate,  rough but sufficient for our analysis, on the measure of discrete cylinders. Considers times \( 0 < t_1 < \cdots < t_n \le 1 \)  and Borel sets \( A_1, \dots, A_n \subset \mathbb{R} \) . 
    Define the corresponding  cylinder set by 
    \[
    \mathcal{C} = \left\{  \omega \in \Omega  \;\middle|\; \omega(t_1) \in A_1, \dots, \omega(t_n) \in A_n \right\}.
    \]
    By the finite-dimensional distribution,  the Wiener measure \( W \) of this set $\mathcal{C}$ is  given by
    \begin{equation}\label{gen form cyl}
    W(\mathcal{C}) = \int_{A_1 \times \cdots \times A_n}
    \prod_{j=0}^{n-1} \frac{1}{\sqrt{2\pi (t_{j+1} - t_{j})}} 
    \exp\left( -\frac{(x_{j+1} - x_{j})^2}{2(t_{j+1} - t_{j})} \right) dx_1 \cdots dx_n,
    \end{equation}
    with the convention \( x_0 = 0 \) and \( t_0 = 0 \). 
    
    Actually we will focus on the following special cylinders given by a path  $ \bar{\omega } \in \Omega$, an integer $ n \ge 1 $ and a real positive number $ \varepsilon > 0$: 
    \begin{equation}
    \mathcal{C}_n ( \bar{\omega} , \varepsilon) := \left\{ \omega   \in \Omega :   \left\vert \omega \left( \frac{j}{n} \right) -  \bar{ \omega} \left( \frac{j}{n} \right) \right\vert \le \varepsilon ,  \text{ for } 1 \le j \le n  \right\} \;. 
    \end{equation}
    This set $\mathcal{C}_n(\bar\omega, \varepsilon)$ will be called the cylinder around $ \bar{\omega}$ of depth $n$ and radius $\varepsilon$.
    Using the formula \eqref{gen form cyl} for $ \mathcal{C}=\mathcal{C}_n ( \bar{\omega} , \varepsilon)$, making the change of variables $ s_i = x_{i} + \bar{\omega} \left( \frac{i}{n} \right) $ for $ i = 1 , \dots , n $, and using the notations $$ \nabla_i s  :=s_{i} - s_{i-1} \qand   \nabla_i \bar{\omega} := \bar{\omega}\left( \frac{i}{n} \right)-  \bar{\omega} \left( \frac{i-1}{n} \right),  $$
    we obtain 
    \begin{equation} \label{eq:cyl-markov}
       W (\mathcal{C}_n ( \bar{\omega} , \varepsilon)) = \left(  \frac{n}{2\pi} \right)^{n/2}   \int_{ [-\varepsilon,   + \varepsilon]^n}  \exp \left( - \frac{n}{2}\sum_{i=1}^n \left[   \nabla s_i - \nabla_i \bar{\omega}\right]^2 \right) \, d s_1 \dots d  s_n  \; , 
    \end{equation}
    with the convention $s_0 = 0 $.

    Observe that $B ( \bar{\omega} , \varepsilon ) \subset \mathcal{C}_n(\bar{\omega} , \varepsilon  )  $. So, we can use the estimation of cylinders' measures to estimate the measures of small balls. 
    
    \begin{lemma}\label{lem:cylinder}
    Let $ (\varepsilon_n)_{ n \ge 1} $ be a positive sequence decreasing to $0$ such that $ \log \varepsilon_n \sim \log \varepsilon_{n+1}$. Then for every $ \bar{\omega} \in \Omega$ we have 
    \begin{equation*}
    \underline{\ord}_{\loc}  W ( \bar{\omega} ) \ge  \liminf_{ n \to + \infty} \frac{\log ( - \log W ( \mathcal{C}_{n} (\bar{\omega} ,\varepsilon_n )  )) }{  -  \log \varepsilon_n } 
    \end{equation*}
    and
    \begin{equation*} 
    \overline{\ord}_{\loc}  W ( \bar{\omega} )  \ge  \limsup_{ n \to + \infty}  \frac{\log ( - \log W ( \mathcal{C}_{n} (\bar{\omega} ,\varepsilon_n )  )) }{  -  \log \varepsilon_n }  \; . 
    \end{equation*}
    \end{lemma}
    \begin{proof}
    This is just because  $ W(B ( \bar{\omega} , \varepsilon)) \le W ( \mathcal{C}_n ( \bar{\omega} , \varepsilon) ) $ for every $ n \ge 1 $ and $ \varepsilon > 0 $. Indeed,  for every sufficiently small $ \varepsilon > 0$ there exists an integer $n$ such that $ \varepsilon_{ n+1} < \varepsilon \le \varepsilon_n$. Then
    \begin{equation*}
     \frac{\log ( - \log W (B (\bar{\omega} ,\varepsilon )  )) }{  -  \log \varepsilon }  \ge 
     \frac{\log \varepsilon_{n}}{\log \varepsilon_{n+1}}\cdot
     \frac{\log ( - \log W ( \mathcal{C}_{n} (\bar{\omega} ,\varepsilon_n )  )) }{  -  \log \varepsilon_{n} } . 
    \end{equation*}
    As $ \log \varepsilon_n \sim \log \varepsilon_{ n+1}$ as $ n \to + \infty$, the desired results follows from \eqref{eq:ord} in Remark \ref{rem:ex}.
    \end{proof}
    The reverse inequalities in Lemma \ref{lem:cylinder} are not true in general as the balls are significantly smaller than the corresponding cylinders. To controll the asymptotic behavior of $W ( B ( \bar{ \omega} , \varepsilon) )   $ with $ W ( \mathcal{C}_n ( \bar{\omega} , \varepsilon )  ) $ we shall control discrete local variations of the path $ \bar{\omega }$  with its deviation from its chord. Let $I_j=[(j-1)/n,j/n]$ and let $ \bar L_j$ be the affine interpolation of $\bar \omega$ between the endpoints of $I_j$. Then we define the \emph{maximal deviation from its chord} as: 
    \begin{equation}
    \Delta_n(\bar \omega):=\max_{1\le j\le n}\; \sup_{t\in I_j} \big| \bar \omega(t)- \bar L_j(t)\big| .
    \end{equation}
    As $ \bar{\omega} $ is continuous on $[0,1]$, thus uniformly continuous, it trivially holds the following fact. 
    \begin{fact}
    For every $ \bar{\omega} \in \Omega$: 
     \[
    \Delta_n(\bar \omega) \xrightarrow[n\to\infty]{} 0.
    \]   
    \end{fact}
    The upper bound $W(\bar\omega, \varepsilon)\le W(C_n(\bar\omega, \varepsilon))$ used in Lemma \ref{lem:cylinder} is trivial. The following theorem  provides a key lower bound of the measure of small balls by the measures  of cylinders and the decay rate of the modulus of continuity $ \Delta_n ( \bar{\omega} ) $. 
 In the  proof, we will approximate the trajectories of the Brownian motion by piecewise affine paths. Then Brownian bridges emerge and the  Kolmogorov–Smirnov law for Brownian bridges will be used. 
    \begin{theorem} \label{thm:cyl-ball}
        Let $ \bar{\omega} \in \Omega$. For every $ \varepsilon > 0 $ and every $n$ such that $ \Delta_n ( \bar{\omega} ) < \varepsilon $, we have
        \begin{equation*}
                W( B ( \bar{\omega} , 3 \varepsilon ) ) \ge \exp ( -   2 n^{2} \varepsilon^2 ) \cdot W ( \mathcal{C}_n ( \bar{\omega} , \varepsilon) )  \; . 
        \end{equation*} 
    \end{theorem}
    \begin{proof}
        For $\omega\in\Omega$, write $L_j$ for the affine interpolation of $\omega$ on the endpoints of $I_j$. Then the process: 
    \[
    X^{(j)}(t):=\omega(t)-L_j(t),\qquad t\in I_j,
    \]
    is a Brownian bridge on $I_j$ which, conditioned on the $\sigma$-field
    \(
    \mathcal F_n:=\sigma\left(\omega(j/n):0\le j\le n\right),
    \)
    is \emph{independent across $j$} and has the standard bridge law on length $1/n$. Indeed, this is the standard Markov and Gaussian regression property of Brownian motion: on each interval $I_j$, conditioned on its two endpoints, the path is the sum of its chord and an independent bridge pinned at $0$ at both ends; independence across disjoint intervals holds given $\mathcal F_n$. Also each $ X^{(j)} $ is independent from the whole grid $ \mathcal{F}_n$.
    We then evaluate the deviation of the sample path $ \omega $ from the trajectory $ \bar{\omega} $ for each $t \in I_j $ as: 
    \[
    \omega(t)-\bar\omega(t)
    =\underbrace{X^{(j)}(t)}_{\text{bridge}} +\underbrace{\left(L_j(t)-\bar L_j(t)\right)}_{\text{grid mismatch}}
    +\underbrace{\left(\bar L_j(t)-\bar\omega(t)\right)}_{\text{template curvature}}.
    \]
    Taking suprema over $t\in I_j$ yields the elementary bound
    \begin{equation}\label{eq:sup-triangle}
     \sup_{t\in I_j}|\omega(t)-\bar\omega(t)|
    \;\le\;  \sup_{t\in I_j}|X^{(j)}(t)| \;+\;  \sup_{t\in I_j}|L_j(t)-\bar L_j(t)|
    \;+\; \Delta_n(\bar\omega).
    \end{equation}
    Observe that the chord difference obeys
    \[
    \sup_{t\in I_j}|L_j(t)-\bar L_j(t)|
    \; = \; \max\Big(|\omega((j-1)/n)-\bar\omega((j-1)/n)|,\;|\omega(j/n)-\bar\omega(j/n)|\Big)
    \]
    because the difference of two affine functions is affine and its maximum on the segment is attained at an endpoint , so that \begin{equation} \label{eq:incl-cyl}
    \bigcap_{ 1 \le j \le n }\left\{  \sup_{t\in I_j}|L_j(t)-\bar L_j(t)| < \varepsilon  \right\}  \; = \; \mathcal{C}_n ( \bar{\omega} , \varepsilon ) \;. 
    \end{equation}
    Combining  \eqref{eq:incl-cyl} and \eqref{eq:sup-triangle}  and using $ \Delta_n(\bar\omega) < \varepsilon $ provides: 
    \begin{equation*}
        B( \bar{\omega} , 3 \varepsilon ) \supset \bigcap_{j =1}^n \left\{  \sup_{t\in I_j}|X^{(j)}(t)|   < \varepsilon \right\} \cap \mathcal{C}_n( \bar{\omega} , \varepsilon   )  
    \end{equation*}
    Using the conditional independance of the bridges conditionnally to $ \mathcal{C}_n ( \bar{\omega} , \varepsilon ) $ and then by independance of each $X^{(j)} $ to this grid allow to obtain:
    \begin{align*}
        W \left( B ( \bar{\omega} ,  3 \varepsilon ) \right) &\ge   \prod_{ j =1}^n W \left(  \sup_{t\in I_j}|X^{(j)}(t)|   < \varepsilon  \  \big\vert  \  \mathcal{F}_n  \right)   \cdot   W \left(  \mathcal{C}_n( \bar{\omega} , \varepsilon   )  \right)     \\
        &=  \prod_{ j =1}^n W  \left(  \sup_{t\in I_j}|X^{(j)}(t)|   < \varepsilon \right)  \cdot   W \left(  \mathcal{C}_n( \bar{\omega} , \varepsilon   )  \right)  \;.  
    \end{align*}
    Thus we have shown: 
    \begin{equation} \label{ineq:cyl-ball}
       W \left( B ( \bar{\omega} ,  3 \varepsilon ) \right) \ge W  \left(  \sup_{ 0 < t \le h }|X(t)|   < \varepsilon \right)^{n}  \cdot   W \left(  \mathcal{C}_n( \bar{\omega} , \varepsilon   )  \right)
    \end{equation}
         where $X$ is a standard Brownian bridge on $[0,h]$ (i.e.\ $X_0=X_h=0$), recalling $h=1/n$. It remains to compute the tail estimate of Brownian bridge. The following series identity is classic (Kolmogorov--Smirnov law for bridges  \cite{an1933sulla}. Also see \cite{Doob1949}
         and \cite{GLS2019}): 
    \[
    W \left(\sup_{t\in[0,h]} |X_t|\ge  \varepsilon \right)
    = 2\sum_{k=1}^\infty (-1)^{k-1} e^{-2k^2 \varepsilon^2/h}
    \; . 
    \]
     As it is alternating with decreasing positive terms, the sum is bounded from below by  $  2 (  e^{- 2 \varepsilon^2/h}  - e^{- 8 \varepsilon^2/h}  )   $ which is at least $e^{-2 \varepsilon ^2/h}  $ for $h$ small enough. Implementing this minoration into \eqref{ineq:cyl-ball} leads to the desired result: 
     \begin{equation*}  
       W \left( B ( \bar{\omega} ,  3 \varepsilon ) \right) \ge  \exp ( - 2 \varepsilon^2 n^2 )  \cdot   W \left(  \mathcal{C}_n( \bar{\omega} , \varepsilon   )  \right)
    \end{equation*}
    \end{proof}

    \subsection{Upper bound of the local exponent $ \overline{\ord}_\loc  W ( \bar{\omega} )$  using the Hölder exponent of $\bar\omega$} 
    Now, having Theorem \ref{thm:cyl-ball}, we can provide the following majoration of the local order of the Wiener measure for non-typical Hölder trajectories.

    \begin{theorem} \label{thm:ubound}
    Let $ \alpha \in (0,1/2]$. Then for every $ \bar{\omega} \in C_\alpha^-$, it holds: 
     \[   \overline{\ord}_\loc  W ( \bar{\omega} ) \le 2 \left(  \alpha^{-1}  - 1 \right) \; . \]  
    \end{theorem}
    
    \begin{proof}
    Fix $ 0 < \beta < \alpha $. 
    Fix $n$.  The formula  \eqref{eq:cyl-markov} with $\varepsilon := n^{-\beta} $ provides: 
    \begin{equation} \label{expr cyl}
    W ( \mathcal{C}_n ( \bar{\omega} , n^{-\beta} )) = \left(  \frac{n}{2\pi} \right)^{n/2}   \int_{ [-n^{-\beta} , + n ^{-\beta} ]^n}  \exp \left( - \frac{n}{2}   \sum_{i=1}^{n} \left[ 
     \nabla s_i - \nabla \bar{\omega}_i \right]^2 \right)  \, d s_1 \dots d  s_n, 
    \end{equation}
    with $s_0 =0$. As $ \bar{\omega} $ lies in $ C^-_{\alpha} $, for large values of $ n$ and every points $ x,y \in [0,1]$ with $ \vert x-y \vert \le n^{-1} $ the increments are controlled by: 
    \begin{equation}
    \vert	 \bar{\omega} ( x) - \bar{\omega} (y) \vert \le n^{-\beta} \;. 
    \end{equation}
    It translates as $ \Delta_n ( \bar{\omega}) \le n^{-\beta}$, thus by Theorem \ref{thm:cyl-ball}, with $ \varepsilon = n^{-\beta} $, it holds: 
    \begin{equation} \label{eq:ball-lbound}
        W ( B ( \bar{\omega} , 3  n^{-\beta} ))  \ge \exp \left( - 2  n^{2 ( 1 - \beta)  }   \right) \cdot W \left( \mathcal{C}_n ( \bar{\omega} , n^{-\beta} )  \right) . 
    \end{equation}
    We shall then estimate $W \left( \mathcal{C}_n ( \bar{\omega} , n^{-\beta} ) \right) $. 
    Note that for every  $i \in \left\{ 1 , \dots, n\right\}$ and for  $s_{i}, s_{i-1} \in [ - n^{-\beta} ,   n ^{-\beta} ] $, it holds:  
    \begin{equation}
        \left\vert	 \nabla s_i - \nabla \bar{\omega}_i  \right\vert   
        \le \vert	s_{i}  \vert + \vert  s_{i-1}  \vert  + \left\vert \nabla \bar{\omega}_i  \right\vert \le 3 n^{-\beta} \; .
    \end{equation}
    By \eqref{expr cyl}, it follows: 
    \begin{align*}
    W \left( \mathcal{C}_n \left( \bar{\omega} ,  n^{-\beta} \right) \right)   &\le \left(  \frac{n}{2\pi} \right)^{n/2}   \int_{ [-n^{-\beta} , + n ^{-\beta} ]^n}  \exp \left( - \frac{n}{2}\sum_{i=1}^{n} \left( 3 n^{-\beta}  \right)^2 \right) \, d s_1 \dots d  s_n  \\
    &=\left( \sqrt{\frac{2}{\pi}} \cdot  n^{1/2-\beta}\right)^{n} \cdot \exp \left( - \frac{9}{2} n^{ 2 ( 1 - \beta) }   \right)   \; . 
    \end{align*}
    Using the above inequality and \eqref{eq:ball-lbound} provides easily: 
    \[  \overline{\ord}_{ \loc} W ( \bar{\omega} )  = \limsup_{ n \to + \infty} \frac{ \log - \log W ( B ( \bar{\omega} ,  3 n^{-\beta} ))   }{  \beta \log n}  \le   2 ( \beta^{-1} - 1 )  \; ,  \]
    where the first equality comes from the sequential characterization of scales  (\cite{helfter2025scales}, Lemma $2.3$). Or see the proof of Lemma \ref{lem:cylinder} for a smiliar argument. Taking $ \beta$ close to $\alpha$ concludes the proof. 
    \end{proof}    
    \subsection{Lower bound of the local order \(\underline{\ord}_{\loc} W( \bar{\omega})\) using Orey exponent of $\bar\omega$}
   We provide here an explicit lower bound of the decay rate of small balls using  Orey exponents of their centers.  
    We first provide the following universal lower bound. 
    \begin{proposition} \label{prop:univ}
    The local orders of $W$ at any trajectory $ \bar{\omega} \in \Omega$ are least $ 2$, namely
    \[     \underline{\ord}_{\loc} W( \bar{\omega}) \ge 2 \; . \] 
    \end{proposition}
    \begin{proof}
    Let us denote $\mathcal{C}_n :=  \mathcal{C}_n ( \bar{\omega} , n^{-1/2})  $. In order to prove that $ \underline \ord_{\loc} W (\bar{\omega}) \ge 2 $, by  Lemma \ref{lem:cylinder}, it suffices to show: 
    \[ \liminf_{ n \to + \infty} \frac{ \log - \log  W ( \mathcal{C}_n ) }{ -  \log n^{-1/2} }  \ge   2 \; . \]
    
    The distribution of the Wiener measure on the cylinder $\mathcal{C}_n$ is   provided by \eqref{eq:cyl-markov} as: 
    \[ W ( \mathcal{C}_n) = \left(  \frac{n}{2\pi} \right)^{n/2}   \int_{ [-n^{-1/2},  n^{-1/2}]^n}  \exp \left( - \frac{n}{2}\sum_{i=1}^n \left[ \nabla s_i - \nabla \bar{\omega}_i \right]^2 \right) \, d s_1 \dots d  s_n . \] 
 Since the term in the exponential is non positive, it follows:
    \[ W ( \mathcal{C}_n ) \leq (2 n^{-1/2} )^n \cdot\left(  \frac{n}{2\pi} \right)^{n/2} =\left( \frac{2  }{ \pi}\right)^{n/2} .\]
    
    Thus we have:
    \[ \liminf_{ n \to + \infty } \frac{ \log - \log  W ( \mathcal{C}_n ) }{ - \log n^{-1/2} }  \ge  \lim_{ n \to + \infty } \frac{ \log n - \log2 + \log \log (\pi / 2 )  }{\frac{1}{2} \log n }  = 2 \; . \] 
    \end{proof}

     The following lemma provides a more precise upper bound of the cylinder measure $W \left( \mathcal{C}_n ( \bar{\omega} , \varepsilon ) \right)$ using the quadratic variation $Q_n(\bar{\omega})$ of the path $\bar{\omega}$.
    
    \begin{lemma} \label{lem:quad-asymp}
    Let $ \bar{\omega} \in \Omega$ be a path. For every integer $ n \ge 1$ such that $ Q_n ( \bar{\omega} ) > 0$ and for every $ \varepsilon > 0$,  the mass of the cylinder  $ \mathcal{C}_n ( \bar{\omega} , \varepsilon) $ verifies:
    \[  - \log W ( \mathcal{C}_n ( \bar{\omega} , \varepsilon )) \ge    n^2 Q_n^2 ( \bar{\omega}) \left( \frac{1}{2} - \frac{2 \varepsilon  }{ Q_n( \bar{\omega})} \right)  -   \frac{n}{2} \log \left(\frac{2n\varepsilon^2 }{\pi  } \right)    \; .   \]
    \end{lemma} 
    
    \begin{proof}
Recall that \eqref{eq:cyl-markov}  provides
    \begin{equation} 
      W (\mathcal{C}_n ( \bar{\omega} , \varepsilon)) = \left(  \frac{n}{2\pi} \right)^{n/2}   \int_{ [-\varepsilon,   + \varepsilon]^n}  \exp \left( - \frac{n}{2}\sum_{i=1}^n \left[   \nabla s_i - \nabla \bar{\omega}_i \right]^2 \right) \, d s_1 \dots d  s_n  \; .  \end{equation} 
    Thus we develop 
    \begin{align*}
     &W (\mathcal{C}_n ( \bar{\omega} , \varepsilon)) \\
      &\le \left(  \frac{n}{2\pi} \right)^{n/2}   \int_{ [-\varepsilon,   + \varepsilon]^n}  \exp \left( - \frac{n}{2} \sum_{i=1}^n \nabla \bar{\omega}_i^2    +  n \sum_{i=1}^{n}\left\vert\nabla s_i  \right\vert \cdot \left\vert \nabla \bar{\omega}_i \right\vert  \right) \, d s_1 \dots d  s_n \\ 
    &\le \left(  \frac{n}{2\pi} \right)^{n/2}   (2\varepsilon)^n \exp \left( - \frac{1}{2} n^2  Q_n^2 ( \bar{\omega} )   +  2 \varepsilon n   \sum_{i=1}^{n} \left\vert
    \nabla \bar{\omega}_i \right\vert  \right) \; .   
    \end{align*}
    Now by Cauchy-Schwartz inequality, it holds: 
    \[  \sum_{i=1}^{n} \left\vert \nabla \bar{\omega}_i \right\vert \le  \sqrt{ n  \sum_{i=1}^{n}  \nabla \bar{\omega}_i^2}  =    n Q_n ( \bar{\omega}) \; . \] 
    Thus combining all the above inequalities together, provides: 
    \[  W (\mathcal{C}_n ( \bar{\omega} , \varepsilon ) )  \le  \left(  \frac{2n \varepsilon^2}{\pi} \right)^{n/2} \cdot \exp \left( - n^2 Q_n^2 ( \bar{\omega})   \left( \frac{1}{2} - \frac{2 \varepsilon}{Q_n (\bar{\omega} ) } \right)  \right) \; .\]
    To conclude the proof it suffices to  apply $- \log$ on both sides of the latter inequality.   
    \end{proof}
    Using the above Lemma \ref{lem:quad-asymp} in which the quadratic variation is involved,
    we obtain the following lower bound of local order of Wiener measure  in terms of the Orey exponents of the trajectory. 
    \begin{theorem} \label{thm:lbound}
    For every $ \bar{\omega} \in \Omega$, we have 
    \[ \underline{\ord}_{\loc} W ( \bar{\omega} ) \ge 2 \frac{1-q^+ ( \bar{\omega}) }{ q^- ( \bar{\omega} ) } \in [0, +\infty ] \; . \]  
    In particular, if $ \bar{\omega} \in \Omega $ admits its Orey exponent $ q(\bar\omega)=\alpha > 0 $, then
    \[  \underline{\ord}_{\loc} W ( \bar{\omega}) \ge 2 ( \alpha^{-1} - 1 ) \;. \]
    \end{theorem}
    
    \begin{proof}
    If $ q^+( \bar{\omega}) = 0 $ or $ q^-( \bar{\omega})  \ge  1/2 $, then the inequality trivially holds, eventually by Proposition \ref{prop:univ}. Thus we can consider real positive numbers $ \alpha, \beta $ such that  $  0 <  \alpha < q^-( \bar{\omega}) \le q^+( \bar{\omega})  < \beta $. For sufficiently large integer $n$, the value of $ Q_n ( \bar{\omega})$ is controlled by:
    \begin{equation}
    n^{-\beta} < Q_n ( \bar{\omega}) < n^{-\alpha} \; . 
    \end{equation}
    Let $ \varepsilon_n := \frac{1}{8} n^{-\beta} $, so that $ \frac{2 \varepsilon_n}{Q_n( \bar{\omega} )} \le 1/4$. 
    Applying Lemma \ref{lem:quad-asymp} to $ \varepsilon = \varepsilon_n$  provides: 
    \[  - \log W ( \mathcal{C}_n ( \bar{\omega} , \varepsilon_n )) \ge   \frac{n^{2 - 2\alpha}}{4}   + O ( n \log ( n))    \;    \]
    as $ n \to + \infty$. 
    It then follows:  
     \begin{equation}
      \liminf_{ n \to + \infty} \frac{\log \left( - \log  W ( \mathcal{C}_n ( \bar{\omega} , \varepsilon_n ) ) \right)}{ - \log \varepsilon_n	 }  \ge \frac{2(1 - \alpha) }{\beta}  \; . 
     \end{equation}
    By Lemma \ref{lem:cylinder},  taking $ \alpha $ close to $ q^{-} ( \bar{ \omega } ) $ and $ \beta $ close to  $ q^+ ( \bar{ \omega}) $ provides the desired result. 
    \end{proof}
    
    \subsection{Proofs of  Theorems \ref{thm:mf-spec},\ref{thm:fBm-loc-ord} and \ref{thm:main}}\label{lieu:proofA}
    After the previous preparations  in this section, 
    we shall conclude the computation of the local order of the Wiener measure (Theorem \ref{thm:main}), its relationship with fBm (Theorem \ref{thm:fBm-loc-ord}) and the evaluation of its multifractal spectrum (Theorem \ref{thm:mf-spec}). We first evaluate the Hausdorff order of a fBm at a typical point. This follows directly from two known results. 
First, the asymptotics for the small-ball probabilities of centered fBm are given in \cite{werner1998existence} and also in \cite{shao1993note, monrad1995small}. Secondly, a result of Dereich--Lifshits \cite{dereichproba05} allows one to extend these asymptotics to balls centered at almost every point.
    \begin{proposition} \label{prop:ord-loc-fBm}
    Let $ B^H $ be a fractional Brownian motion with  Hurst parameter $ H \in (0 ,1) $, then the local order of the measure $W^H $ for the uniform norm is 
    \[ \ord_{\loc} W^H(\bar{\omega})  = \frac{1}{H}, \quad  W^H\text{-a.s.} \;. \] 
    \end{proposition}
    \begin{proof}
    By the main result due to Li and Werner in \cite{werner1998existence}
    (weaker results of Monrad and Rootz\'en \cite{monrad1995small} or of Shao \cite{shao1993note} will do too), there exists a constant   
    $\kappa_H \in (0, \infty)$ such that:
    \begin{equation} \label{small ball fBm}
       \quad   \lim_{ \varepsilon \to 0 } - \varepsilon^{1/H} \log W^H ( B(0,\varepsilon))  =  \kappa_H \;. 
    \end{equation}
    Then, by Theorem $2.1$ and Remark $2.2$ in \cite{dereichproba05}  it holds: 
    \begin{equation} \label{dl det}
    \text{a.s.}  \quad 
     - \log W^H ( B( \bar{\omega} ,\varepsilon))   \le  -  \log W^H ( B(0,\varepsilon))   \lesssim  - 2  \log W^H ( B( \bar{\omega} , \varepsilon/2))  \quad \text{as $ \varepsilon \to 0$} . 
    \end{equation}
    Then the claimed result follows from \eqref{dl det} and \eqref{small ball fBm}.
    \end{proof}
    As we have pointed out above, weaker results than that of \eqref{small ball fBm}, but  sufficient for our purpose, have been previously obtained in \cite{shao1993note,monrad1995small}.
    These weaker results can be similarly proven for $ L^p$-norms $ p > 1$ in place of  the uniform norm.  
    
    \bigskip
    {\bf Proof of Theorems \ref{thm:main} $\&$ \ref{thm:fBm-loc-ord} \label{lieu:proofB}.
    }
    We prove simultaneously both theorems.
    Let  $ \bar{ \omega } \in C_\alpha^\crit $ ($0<\alpha <1$). Then by Theorem \ref{thm:lbound}: 
    \[ \underline{\ord}_\loc W ( \bar{\omega} )  \ge 2 ( \alpha^{-1} - 1) \;. \]
    Conversely by Theorem \ref{thm:ubound}:
    \[ \overline{\ord}_\loc W ( \bar{\omega} )  \le 2 ( \alpha^{-1} - 1) \;. \]
    We have thus established the first assertion in Theorem \ref{thm:main}.
    
    \smallskip
    According to Proposition \ref{prop:reg-fBm},  $ B^H $ lies in $ C_H^\crit $ almost surely.  We obtain Theorem \ref{thm:fBm-loc-ord}. 

    \smallskip
    Now,  by  Corollary \ref{cor:spec} and Proposition \ref{prop:ord-loc-fBm},  we first get
    \begin{equation}
         \ord_{H} C_H^\crit \ge \ord_{H} W^H  = \frac{1}{H} \;.  
    \end{equation}
    For the reverse inequality, note that  $ C_H^- \subset C_\alpha$ if $ \alpha < H $. By Theorem E in \cite{helfter2025scales},  the Hausdorff order for the uniform norm of the unit ball in $C^\alpha$ is equal to $\alpha^{-1}$.  As $C^\alpha $ is separable for the uniform norm it is covered by a countable union of shifts of the unit ball, and as Hausdorff order is $\sigma$-stable, it follows that indeed $ \ord_H C^\alpha = \frac{1}{\alpha}$. Consequently: 
    \begin{equation}
        \ord_H C_H^\crit \ge  \frac{1}{\alpha} \quad \forall \alpha < H\;.
    \end{equation}
    Thus we indeed obtain $ \ord_H C_H^\crit = \frac{1}{H}$. 
    Then we conclude the proof of Theorem \ref{thm:main} by inferring the result of Proposition \ref{prop:univ}. \hfill $\square$

\bigskip
    {\bf Proof of Theorem \ref{thm:mf-spec}.}
    First note that the second assertion of Theorem \ref{thm:mf-spec} follows directly from the first one. For this first assertion, as $ \Omega^\crit = \bigcup_{ \alpha > 0} C_\alpha^\crit $ and as Theorem \ref{thm:main} yields $ \ord_{\loc} W ( \bar{\omega} ) =  2 ( \alpha^{-1} - 1)  $  while $ \ord_H  C_\alpha^\crit  $ for $ \alpha  \le 1/2 $. It follows that for $\xi = 2 ( \alpha^{-1} - 1) $,  the level set $E^\crit(\xi) $ in  $\Omega^\crit $ coincides with $ C_\alpha^\crit$ so that $ \ord_H E^\crit(\xi)  = \frac{1}{\alpha}= 1 + \frac{\xi}{2} $. This treats the case $ \xi \ge 2$. We conclude the proof by noticing that for $ \xi < 2$ it holds $E(\xi) = \emptyset $ by Proposition \ref{prop:univ}. \hfill $\square$

\section{Questions and open problems} \label{sec:quest}
The multifractal formalism introduced in this work opens several perspectives toward a general theory of multifractal analysis in infinite-dimensional settings. While our present analysis focuses on the Wiener measure as a prototypical example of a measure supported on an infinite-dimensional space, the notions of multifractal spectra at fixed scales  developed here extend beyond the Gaussian framework and should apply to a broad class of processes.  
In this section, we discuss several open problems and conjectures that, in our view, delineate the next natural steps toward a comprehensive multifractal theory for infinite-dimensional measures.

\subsection{Beyond $\Omega^\crit $} \label{lieu:quest-omega-crit}
A first natural question arising from Theorem \ref{thm:mf-spec} is whether the Hausdorff order of a level set on the whole space $ \Omega $ coincides with the Hausdorff order of its restriction to the space $ \Omega^\crit$. 
Since $\Omega $ is not $ \sigma$-compact and that we do not know  whether $ \ord_\alpha$ verifies the increasing sets lemma,  we can not apply Corollary \ref{coro:mass-distrib-pple}. 
However the subspace $  C_0^+ := \bigcup_{ 0 < \alpha < 1 } C_\alpha $ of  functions admitting at least one Hölder exponent is $\sigma$-compact and contains $ \Omega^\crit$. In view of Corollary \ref{coro:mass-distrib-pple-intro}, to show that $ \Omega^\crit$ can be replaced by $ C_0^+ $ in the statement of Theorem \ref{thm:mf-spec},  it would be enough to answer to the following question. 
\begin{question}
  Consider  $ \alpha >  2 $ and let $ X $ be a random vector such that almost surely $X \in  C_0^+$ and $  X \in E_W^{\ord} ( \alpha)  $. 
  Does it hold $$ \ord_H^* \mathcal{\mu} \le 1  + \frac{\alpha}{2}$$
  where $ \mu $ is the probability law of $X$ ? 
\end{question}

Besides, the question of the spectrum of Wiener measures can naturally be asked for scales other than the order. Indeed, $\Omega$ is an infinite-dimensional Banach space for which, in view of \cite{hel25}, there exist subsets of arbitrary scales. One may then ask whether these subsets can be realized as level sets of highly irregular trajectories of the Wiener measure.

\subsection{Multifractal analysis of fBm}
 In view of Proposition \ref{prop:ord-loc-fBm} and our results in Theorems \ref{thm:main} and \ref{thm:fBm-loc-ord}
on the local order of the Wiener measure, combined with the characterization of fBms in terms of their Hölder and Orey exponents in Proposition \ref{prop:reg-fBm},  a broader phenomenon is strongly suggested.  We expect that the local order of a given fractional Brownian measure can  be described using trajectories of other fBm's.
    \begin{question} \label{conj:spec-fBm}
    Is the following conjecture true ?  \medskip
        For every $ H^* \in ( 0,1) $, the local order of the fractional Brownian measure $W^{H^*} $  in the uniform norm evaluated at $B^H$ for $ H \in (0,1) $ is $ W^H$-a.s.
        \begin{equation*}
            \ord_{ \loc  } W^{H^*} ( B^H) =
            \begin{dcases}
                \frac{1+ 2 ( H^{*} - H)}{H} & \qquad  \text{ if } \, H < H^* \\
                \frac{1}{H^*} & \qquad  \text{ otherwise.}
            \end{dcases}
        \end{equation*}
    In particular on $ \Omega^\crit $ it would hold: 
    \begin{equation*}
        f_{W^{H^*}}^\ord ( \xi) = 
        \begin{dcases}
        \frac{\xi  + 2  }{ 1 + 2 H^* } & \text{ if } \xi \ge \frac{1}{H^*} \\
             - \infty   & \text{ otherwise. }  
        \end{dcases}\;. 
    \end{equation*}
    \end{question} 
 This conjecture would unify all fractional Brownian measures under a single multifractal framework, where the interplay between the two indices $H$ and $H^*$ is described by the geometry of local fluctuations of trajectories. 

\smallskip
The main difficulty is that a general fBm has correlated increments, whereas the independence of increments for standard Brownian motion was used intensively in our proofs.
\subsection{Independance of the metric}
    
Our analysis has been carried out only under the uniform norm.  Recall that  many results of  Dereich-Lifshits in \cite{dereichproba05} are valid for a wide range of norms such as $ L^p$-norms for $ p >1$. Consequently, it seems  natural to ask whether the result Theorem \ref{thm:mf-spec} holds true when the uniform norm is replaced by $ L^p $-norm. 
 \begin{question}
Is it true that the local order and multifractal spectrum of $W$ are independent of the choice of $p > 1$, i.e.
\begin{equation}
 \ord_{\loc} (W , \| \cdot \|_p ) = \ord_{\loc} (W , \| \cdot \|_\infty )
\end{equation}
and,  in particular,  that the order spectrum is universal across $L^p$-norms on $\Omega$ ? 
\end{question}
A positive answer would indicate that the geometry of infinite-dimensional measures such a Wiener measure is not so sensitive to the underlying functional structure, perhaps revealing that anisotropies in Gaussian measures is somehow tied to topological properties. 
More generally we can ask Question \ref{conj:spec-fBm} for any $ L^p$-norm with $ p > 1$.

\subsection{Possible multifractal spectra}
As mentioned earlier, the multifractal formalism developed here naturally extends to other (Gaussian) processes. A  notable problem
that arises, is to determine whether different covariance structures lead to distinct scaling behaviors. Indeed, scales provide a natural framework to encompass the study of other Gaussian processes, such as Ornstein–Uhlenbeck, Bessel, or Lévy-driven processes. A fundamental problem is to determine whether distinct covariance structures could lead to distinct scaling behaviors.
\begin{question}\label{quest:gaussian}
Given a scaling $ \Scl $, does there exist a centered Gaussian process  on $[0,1]$ with continuous sample paths such that the multifractal spectrum of its law at scale $ \Scl $ is not trivial ? 
\end{question}

    \subsection{Relationship between local scales and dimensions of paths}
    There is a well-known connection between scaling exponents and the fractal dimension of the graph. Let us recall the following well known result.
    See \cite{Kahane1985}
    (Theorem 6 in Chap. 10, p. 139 and    Theorem 7 in Chap. 18, p.278). 
    Let $f:[0,1]\to\mathbb{R}$ be $\alpha$-Hölder continuous with $0<\alpha<1$. Then
    \[
    \dim_H \mathrm{Graph}(f)\;\le\;\overline{\dim}_B \mathrm{Graph}(f)\;\le\;2-\alpha.
    \]
    Moreover, for a typical trajectory of fractional Brownian motion $B^H$ with Hurst index $H$, this upper bound is sharp, and a.s. we have
    \[
    \dim_H \mathrm{Graph}(B^H)=\dim_B \mathrm{Graph}(B^H)=2-H.
    \]
    
    \begin{question}\label{quest:graph}
Given a Gaussian process with law $\mu$, does there exists a function $F$ such that for every other Gaussian process $X$ with continuous sample paths,
\[
\ord_\loc \mu(X) = F\!\left(\dim_H(\mathrm{Graph}(X))\right) \qquad \text{ a.e.}  \;  ? 
\]
\end{question}
Such a result would establish a geometric link between the multifractal order of the probability law and the Hausdorff dimension of its realizations.
Note that, in particular, for the standard Brownian motion the only candidate to answer Question \ref{quest:graph} is: 
\[
F  : x \in [1,2) \mapsto  \frac{2x - 2 }{2 -x } \;.
\]

Here we refer to Theorem \ref{thm:fBm-loc-ord} with $ H = 2 -x $.  
\subsection{Thermodynamical formalism}
Finally, we must mention that the generalization of \emph{coarse multifractal spectra} (see for instance Chapter~17 of~\cite{falconer2004fractal}) to arbitrary scales is not addressed in the present work. 

In the dimensional case, to estimate \( f(\alpha) = \dim_H E(\alpha)  \), one introduces the \emph{partition function}. Let \( \mathcal{D}_j \) be the collection of dyadic cubes of side length \( 2^{-j} \) in \( [0,1]^d \). For \( q \in \mathbb{R} \), define:
    \[
    Z(q,j) := \sum_{Q \in \mathcal{D}_j, \mu(Q)\not=0} \mu(Q)^q,
    \]
    and the associated $\tau$-function :
    \begin{equation}  \label{eq:legendre}
    \tau(q) := \liminf_{j \to \infty} \frac{\log Z(q,j)}{\log 2^{-j}}.
    \end{equation}
    The multifractal formalism, suggested by physicists in the context of turbulence~\cite{benzi1984multifractal,frisch1985turbulence}, asserts that the multifractal spectrum \( f(\alpha) \) is given by the Legendre transform of \( \tau(q) \), i.e.,
    \begin{equation} \label{eq:multifractal-formalism}
    f(\alpha) \leq \tau^*(\alpha) := \inf_{q \in \mathbb{R}} \left( \alpha q - \tau(q) \right).  
    \end{equation}
    Moreover the equality holds for regular self similar measures.
It appears that extending such definition and finding a multifractal formalism for coarse spectrum is not straightforward at arbitrary  scale, even in the specific case of order. It would be of great interest to develop an analogue of the multifractal formalism, as expressed in~\eqref{eq:multifractal-formalism}, valid at any scale. 
Such an extension could reveal new forms of duality between scales and corresponding generalized spectra.
    \begin{question}
       Can we find an analogue of \eqref{eq:multifractal-formalism} for a general scaling $ \Scl $ ? 
    \end{question}

\bibliographystyle{alpha}
\small
\bibliography{ref.bib}
    
\end{document}